\documentclass[10pt,a4paper,oneside,reqno]{amsart}

\usepackage{amsmath,amsfonts}

\usepackage{multirow}
\usepackage{float}

\usepackage{geometry}
\usepackage{cancel }

\usepackage{graphicx}
\usepackage{pgfplots}
\pgfplotsset{compat=newest}
\usepackage{tikz}
\usetikzlibrary{decorations.markings,shapes.geometric,arrows,calc}

\newtheorem{remark}{Remark}
\newcommand{\dsum}{\displaystyle\sum}

\usepackage{hyperref}

\let\origmaketitle\maketitle
\def\maketitle{
  \begingroup
  \def\uppercasenonmath##1{} 
  \let\MakeUppercase\relax 
  \origmaketitle
  \endgroup
}
\begin{document}

\title[Line Planning and Timetabling in Subway Networks]{\large{An optimization model for line planning and timetabling in automated urban metro subway networks}}

\author{V\'ictor Blanco}
\address{IEMath-GR, Universidad de Granada.}
\email{vblanco@ugr.es}

\author{Eduardo Conde}
\address{Dep. Statistics \& OR, Universidad de Sevilla.}
\email{educon@us.es}

\author{Yolanda Hinojosa}
\address{Dep.  Applied Economics I, Universidad de Sevilla.}
\email{yhinojos@us.es}

\author{Justo Puerto}
\address{Dep. Statistics \& OR, Universidad de Sevilla.}
\email{puerto@us.es}

\thanks{The authors were partially supported by the projects FQM-5849 (Junta de Andaluc\'ia $\backslash$ FEDER),   MTM2016-74983-C2-1-R (MINECO, Spain) and contract 1853/0257 (Soci\'et\'e Metrolab{\textregistered}, Service Contr\^ole de Gestion).}

\date{\today}

\begin{abstract}
In this paper we present  {a} Mixed Integer Nonlinear Programming model that we developed as part of a pilot study requested by the R\&D company {\sc Metrolab\textregistered}\footnote{Soci\'et\'e Metrolab\textregistered, Service Contr\^ole de Gestion, registered on the Paris Trade and Companies Register under the number 532 684 685 RCS Paris, with its registered office at 117/119 Quai de Valmy - 75010 PARIS - FRANCE} in order to design tools for finding solutions for line planning and timetable situations in automated urban metro subway networks. Our model incorporates important factors in public transportation systems from both, a cost-oriented and a passenger-oriented perspective, as time-dependent demands, interchange stations, short-turns and technical features of the trains in use. The incoming flows of passengers are modeled by means of piecewise linear demand functions which are parameterized in terms of arrival rates and bulk arrivals. Decisions about frequencies, train capacities, short-turning and timetables for a given planning horizon are jointly integrated to be optimized in our model.  Finally, a novel Math-Heuristic approach is proposed to solve the problem. The results of extensive computational experiments are reported to show its  applicability and effectiveness to handle real-world subway networks.
\end{abstract}

\maketitle

\section{Introduction}

In this paper, we propose a model for line planning and timetabling on general urban subway transportation systems. This study was originated by a real-world problem proposed by Metrolab{\textregistered}, a French R\&D company, dealing with the line planning and timetabling of trains of existing subway networks. It was a pilot experience to automatize the decision making process, at the tactical and operational level, of a small section of the {Paris subway network.}

The development of flexible tools to control, automatically, a transportation system according to a set of indicators of its service quality may have a considerable impact in its efficiency and usefulness. The quality perception of a public transportation system from a customer's point of view is highly dependent of its reliability, comfortability and effectiveness when comparing with alternative transportation means. If just poor-quality connections are offered or the quality-price relationship does not fulfil the passengers' expectations, they may decide to use alternative transportation means. Therefore, the quality of a public transportation system, from the passengers' point of view, is a key objective in its design and management, besides infrastructure constraints, operational limitations or budget considerations~\cite{goe17}.

The existing literature on line planning is very extensive (see e.g., \cite{des07, goe17, sch12}  and the references therein), including different models which can be classified according to the decisions covered (determination of train routes, frequency setting, or both), infrastructure and operational aspects, objective functions and the way in which passengers' decisions are taken into account in the decision making process. For instance, in \cite{goe17}, line planning models are classified with respect to their  objective functions into models with cost-oriented or with passenger-oriented objective functions.  Our model will consider both points of view in an attempt to find an equilibrium between these conflicting objectives which, as mentioned in \cite{gui08}, is an important challenge in a public transportation system.

However, building an effective  model is much more complex than selecting the appropriate nature of the objective function. The correct delimitations of the considered features in the context of the transportation system or in the set of customers served by this system is a very important {phase in the actual  modeling. Often, an initial line planing must be re-engineered  motivated by changes on costumers' flows induced by changes in \emph{passengers' route choices}, \cite{goe17}.} This gives rise to a \emph{bilevel optimization problem} with a line planning problem on the upper level and  a passenger's route choice problem on the lower level. Moreover, the existence of several decision-makers  is not the only difficulty of this problem. There exist uncertainties in the number of passengers that must be served (demands), in the origin-destination pairs that customers want to go across or in the times needed to cover the network links due, for instance, to machine failures or other incidents. In this situation it seems hard to integrate all these elements in  a suitable optimization model.

Usually, some  simplifications must be assumed in order to obtain operational solutions for realistic situations. The usefulness of the model will be strongly conditioned by these assumptions.  Having a valid model for a wide range of scenarios might be a lofty target but one worth aiming for. The resulting approximation can be seen not only as a model to optimize the existing resources in a given transportation system but also as a \emph{what-if} tool to make  rational decisions. For instance, it can be used to check the implementation of  possible modifications in the structure of the existing network (adding stations, new connections, \ldots) or in the conditions (demand variation, service disruptions, \ldots) under which the transportation system is working at present.

Line planning is only one of the planning process's stages. Indeed, following \cite{des07,mic09,sch17}, the planning process in public transportation includes several phases that usually are sequentially executed in the following order:
\begin{enumerate}
  \item {\bf network design}, where the stations, links and routes of the lines are established,
  \item {\bf line planning}, specifying the frequency and the capacity of the vehicles used in each line (\emph{line concept}, \cite{goe17}),
  \item {\bf timetabling}, {defining } the arrival/departure times and
  \item {\bf scheduling},  in which vehicles and/or crews are planned.
\end{enumerate}
The first phase, namely network design, is done at the strategic level and implies a high cost (see e.g., \cite{bpr11}). Moreover, the remaining steps involve decisions at the tactical and operational levels which are conditioned by that design. Thus, it seems appropriate to assess different designs by means of a \emph{what-if analysis} based on  the efficiency of the system under a given scenario. In that case, after initializing with a \emph{reasonable} design, the procedure could \emph{optimize} the efficiency of the transportation system using tactical/operational decisions. This may reveal possible weaknesses of the current design and, after fixing some of them, will give rise to a new design. The process could be repeated until a compromise transportation design is found.

In addition to the literature on line planning commented above, one may also find a rich literature on timetabling and scheduling (see e. g. \cite{bar14, cap02, cap06, lee09, sun14} for timetabling models and \cite{can11, dar07, tek18} for scheduling models). Timetabling models are usually classified  according to the capacity of the transportation system or the requirements of passengers. Regarding to the scheduling literature, models can be classified into single and multiple depot frameworks~\cite{bun09}. Besides, we can find periodic and time-dependent timetabling and scheduling models and many other variants depending on the considered constraints. Considerable effort has also been made to analyze the practical effects of the different elements or features covered by these models. For instance, as mentioned in \cite{niu15}, periodic timetables easily allow passengers to remember the exact departure times at stations but, in general, they are not fully sensitive to time-varying passenger demands, which could result in long waiting times and reduced service reliability, particularly under irregular over-saturated conditions.

Once the framework of the transportation system has been delimited, the resulting model should be optimized. However, as pointed out in \cite{sch17}, going through the above-mentioned stages of the planning process in a sequential way, often leads just to suboptimal solutions. Recent research (see \cite{gor13,gui08,lopp17,opp18,sch17}) is increasingly oriented towards {\bf integrated planning} in which two or even more of these planning stages are simultaneously addressed. These integrated optimization models are frequently superior to those optimizing sequentially the considered stages, as it has been recognized in the literature (see e.g., \cite{sch17} and the references therein).

In this paper we propose a new integrated model in which line planning and timetabling are simultaneously optimized using a combination of cost- and passenger-oriented objective function. In our model, time-dependency on demands is considered in addition to two other elements that, as far as we know, have only been addressed separately in the literature: short-turns and interchange stations.

\textit{Short-turning} is a tactical decision for which some trains can perform short cycles in order to increase frequency in specific sections of a line. In general, due
to their analytical complexity, the approaches to manage short-turnings in the context of railways planning are based on particular cases and there are no general models that can be applied without modifications to every situation~\cite{can16}. Besides, most of the literature analyzing short-turning in a railway context is limited to a single two-way transit line  \cite{can16, cor11, ort09, wee18}.  Our model incorporates decisions concerning the activation of short-turns simultaneously in several lines of the network.

On the other hand, \textit{interchange} or \textit{transfer stations} are shared by several lines in the system allowing  passengers to change from one to another line. Papers dealing with interchange stations (see e.g. \cite{hei18,kan16,won08}) usually aim for minimizing the total transfer waiting time of passengers by synchronizing train arrival times at transfer stations. We will manage here the effective flows of passengers at the interchange stations and compute the corresponding effects both in the quality of the service and in technical constraints, as those related to the capacities of the trains.

The final aim is to model a real system, inspired by one initially proposed by Metrolab{\textregistered}, using its more relevant features whilst its computational tractability is preserved. This is an important challenge taking into account the combinatorial, stochastic, multilevel and multiobjective nature of the problem. The resulting outcome is a Mixed Integer Second Order Cone Programming model, which can be solved using off-the-shelf optimization solvers, but only for limited sizes. For larger sizes we propose a Math-Heuristic approach in which the system is decoupled into different lines. Afterward, each subsystem is optimized individually but including in the input data the flows of passengers generated after optimizing other lines. The process is repeated like in a \emph{block coordinate descent} procedure (see e.g., \cite{ber99}) until a given stoping rule is verified.

The remainder of this paper is organized as follows. Section \ref{sec:1} deals with a detailed description of the optimization problem and its main elements. In Section \ref{sec:2} we present the Mathematical Programming formulation of the problem. The demand function modelling how the flow of passengers entering into certain stations of a line changes according to the effects of the other lines or due to external factors is detailed in Section \ref{sec:3}. The usefulness of the proposed model is illustrated in Section \ref{ex0} with a case study using real data on a section of the Paris subway provided by Metrolab{\textregistered}. Our Math-Heuristic approach is proposed in Section \ref{sec:4} and the corresponding computational results, including its comparison with the exact method, are reported in Section \ref{sec:5} using several network topologies adapted from the literature.  The paper ends with a section where some conclusions and future extensions of the model are outlined.

\section{Problem description}
\label{sec:1}

Let us assume that the technical features of a public transportation network, which is a part of a complex underground train system, are known. Our goal is to model the problem of how to operate different metro lines on this network according to a set of technical requirements and a given structure of the demand requesting for this service. Suppose that a set of routes for the potential lines (\emph{lines pool}) and a set of interchange stations are specified.  Furthermore,  some other factors involved in the system performance such as passenger flows, set of possible train capacities, maximum number of allowed trips in a given line during the planning horizon, demand fluctuation (e.g., rush/off-peak hours) or stopping time windows, amongst others, are assumed to be also available. The goal is to model such a system using its more relevant features whilst its computational tractability is preserved. In the description of our approach we will consider three main blocks: input data, feasible actions and assessment of a particular solution together with some additional specifications  of our model.

\subsection{Input data}

In addition to the input network, including the topological route map and the interchange stations, the main block of input data corresponds to the passenger flow amongst the considered set of stations. Obviously, the number of passengers awaiting in a specific station for the next train which connects to a given destination is a stochastic process. Furthermore, the stochastic processes corresponding to the set of considered stations are  interdependent  due to the flow relationships amongst the stations which, in addition, may change over time. In our model, given that we deal with a  heavily congested subway line, these stochastic processes will be replaced by average rates of the number of passengers.
The dynamic dependence of these processes on time should be preserved in some way in the optimization model since it is one of the relevant elements in order to obtain realistic operating solutions. We will do it through a function measuring the intensity of the demand.

In our framework, we assume that two main situations affect to the passengers flows: transitions between stations or lines and arrival of external demand functions. The first one refers to the behaviour of the passengers with respect to the mobility pattern. These data fix an assignment for the destination of the passengers catching a train in a given station, what is known in the literature as \emph{line planning with route assignment}, \cite{goe17}. The alternative approach of line planning with route choice seems to be more appropriate just in those transportation systems with high density of connections (with alternative paths between two given locations) and low trip frequency. However, this is not the case in our model since these two features rarely appear in underground train systems.

We will use an Origin-Destination (OD)-matrix per line having as entries the proportions of passengers moving between pairs of stations of the line and also a set of values quantifying the proportion of passengers which want to change from one metro line to another in each one of the transfer stations. These proportions may be considered as estimations of the probabilities with which a passenger moves through the network and could be dependent on the dynamic nature of the transportation system. Following \cite{sch17}, in order to model the passengers' flows,  OD-data gives rise to more realistic applications than those based on traffic loads since the paths followed by passengers depend strongly on the line concept. As observed by several authors (\cite{goe17,sch17}), the optimization models derived from the management of OD-data are often harder to be solved numerically, and thus, approximated ad-hoc algorithms need to be used to deal with problems of realistic sizes.

The second concept, the arrival of external demand functions, refers to the intensity of use of the transportation network. The external demand models the incoming flow of passengers entering to the system from outside during the planning horizon. These functions determine the relative importance, in the overall planning cost, of the stations used to access the system  and it may change depending on time. Also, rush hours at given time periods{, irregular weather conditions,} or the celebration of events at certain places close to stations may provoke increasing or decreasing of the incoming flow of passengers taken into account in our model.

These external demand functions, together with the OD-data corresponding to movements between pairs of stations {and the proportions of transfer passengers}, give rise to a model having time-dependent passenger flows.  Time-dependent flows are an appropriate feature of a realistic approach for traffic planning and, as pointed out in \cite{sch17}, at present, there is not much research literature covering integrated optimization planning models under these conditions.

\subsection{Feasible actions}

In our model, the line concept design is specified by choosing the operating frequencies and the train capacities for each line considered in the lines pool. Furthermore, our line concept design allows \textit{short-turning} in some lines, i.e., the  possibility of activating, for some or all the trains, and at certain time periods,  short cycles, in order to increase the frequency in specific (consecutive) stations suffering from intensive demand.  This situation is typical in lines which connect distant residential areas with the city center or economic centers.

In the integrated planning model, the line concept design is optimized together with their corresponding timetable. Both elements define the feasible solutions of the problem once a set of technical constraints, involving passenger flows and leaving/arrival times, is specified.

The line concept design will be the main source of discrete decision variables for the formulation of our problem. As, for instance, the selection, among a finite set of capacities for the trains. On the contrary, the actual number of trips of a given line will be modeled using a finite set of replicas of continuous variables corresponding to potential timetables. On the basis of these decision variables a number of additional auxiliary variables will be considered in the optimization model in order to control the times in which different events happen at each station during the planning horizon. In our model, unlike most of the approaches deriving optimal timetables in transportation systems, the departure times are not discretized and the period of time elapsing between consecutive train departures is not constrained to be constant. Thus, it provides more flexible decisions as well as timetables sensitive to the changeable conditions in the passengers's flow during the planning horizon, turning out, in general, in an aperiodic timetabling.

Different train speeds are not taken into account in our line concept design because the pilot proposal by Metrolab{\textregistered}  only considered constant and fixed speeds between each pair of consecutive stations. Nevertheless, continuous variables modeling the speed of a train between two consecutive stations could be \textit{easily} added to our model as explained in Remark \ref{remark1}.

\subsection{Assessing a planning solution}

As commented above, due to the multi-objective nature of the problem, one of the most difficult modeling issues is that of assessing a feasible solution (line concept+timetable). Maintenance/operational planning costs are usually easy to handle as a part of the objective function. However, in order to consider also the passenger-oriented nature of the objective,  the cost induced by the quality of the service provided by the system should be included in the objective function. This will be incorporated into the model quantifying the cost of unmet (non-served) demand, that is, passengers who cannot take a train due to lack of capacity. Non-served passengers contribute a given amount, in terms of costs, due to their confidence loss in the transportation system, their balking rate or a combination of these two and some other factors. Obviously, calibrating these costs is not an easy task, but it may be partially handled by managing a finite set of cost estimates and solving the problem for each one of them in order to evaluate the influence of these hard-to-calibrate parameters in the proposed solution.

\subsection{Specifications  of our model}
 In the following we list the assumptions that are imposed to derive a suitable Mathematical Programming formulation for our model.
\begin{description}
\item[\textit{Planning horizon}:] Our model considers a continuous time interval in which all the events must start and decisions occur. The range of this interval depends on the data collection accuracy and with-in-day variability of the demand changes. This planning horizon is fixed a priori but, as mentioned in Remark \ref{remark2}, our model allows to join two consecutive planning horizons by passing data about numbers of passengers and arrival/departure times obtained from an optimal solution on a given planning period as input data for the next one.
  \item[\textit{One direction trips}:] We assume that each line in the pool operates only in one direction, from a given line-header station to the final one. Usual round-trip lines are modeled as two symmetric lines by interchanging the order of the stations. A \emph{trip} consists in making the complete walk along all the ordered stations of a given line, from the head to the final station. Thus, a physical round-trip starting and ending at the same station will be given by two lines sharing the stations but in the opposite order.
 \item[\textit{Short-turns}:]  We consider that specific sections of consecutive stations in certain lines are allowed to be activated in some of the trips.
  \item[\textit{Interchange stations}:]  We  consider that the lines in the lines pool   may share common interchange stations where some of the passengers change of line to go to their final destination.
  \item[\textit{Train capacities}:]  The model assumes that a finite set of admissible capacities for the trains operating a line is given.
  \item[\textit{Safety interval}:] A minimum security time window between consecutive trips in any line is established.
  \item[\textit{Maximum number of trips}:] We assume, w.l.o.g., that the maximum number of possible trips in each line is given beforehand. Note that an upper bound of this maximum number can be obtained taking into account the range of the planning horizon and the safety interval.  In our formulation, we resort to a set of decision variables controlling the time in which different events happen. These variables must be replicated as many times as the maximum number of trips, although only some of those \emph{trip-variables} are activated and then represent actual trips. Hence, the size of the formulation strongly depends on this maximum number. The idea is to consider that the potentially variable number of trips of the line concept is fixed to the maximum number although, in fact, some of them are really \textit{fake trips}. This trick will ease the task of building constraints to ensure non-overlapping events and {the safety interval} between consecutive trips in the stations of a given line.
  \item[\textit{Piecewise linear cumulative incoming demand}:] We model the accumulated number of passengers arriving to each station up to a given instant during the planning horizon by the so-called {\it demand function}. With this function we manage variable arrival rates during the planning horizon and bulk arrivals due to special events, like for instance the end of a football match in a close location, the arrival of passengers coming from another transportation mean or, in the interchange stations, the arrival of passengers coming from another line of the subway network. As proposed by Metrolab{\textregistered} for its pilot experience, we consider that this function is a piecewise linear function of time (further details are given in Section \ref{sec:3}).
\end{description}

Figure \ref{fig:00} illustrates some of the considered features of the networks under study. There, we represent by circles or squares the nodes corresponding to the stations of two subway lines. An interchange station common to two lines is marked with a black square. We have also included a {possible} short-turn in the \textit{horizontal} line covering a set of four consecutive squared stations (drawn as a dashed line in the picture).

 \begin{center}
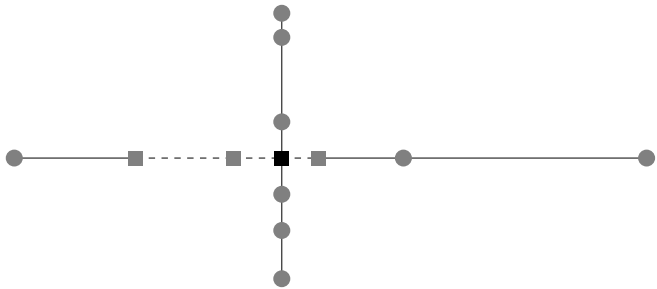
\begin{figure}
\begin{center}
\scalebox{0.8}{\begin{tikzpicture}[scale=0.4]

\coordinate(A) at (0,0);
\coordinate(B) at (5,0);
\coordinate(C) at (9,0);
\coordinate(D) at (11,0);
\coordinate(E) at (12.5,0);
\coordinate(F) at (16,0);
\coordinate(G) at (26,0);
\coordinate(AA) at (11,6);
\coordinate(BB) at (11,5);
\coordinate(CC) at (11,1.5);
\coordinate(DD) at (11,0);
\coordinate(EE) at (11,-1.5);
\coordinate(FF) at (11, -3);
\coordinate(GG) at (11,-5);

\draw (AA)--(BB);
\draw (BB)--(CC);
\draw (CC)--(DD);
\draw (DD)--(EE);
\draw (EE)--(FF);
\draw (FF)--(GG);

\draw (A)--(B);
\draw[dashed] (B)--(C);
\draw[dashed]  (C)--(D);
\draw[dashed]  (D)--(E);
\draw (E)--(F);
\draw (F)--(G);

\fill[gray](AA) circle (0.35);
\fill[gray](BB) circle (0.35);
\fill[gray](CC) circle (0.35);
\fill[gray](EE) circle (0.35);
\fill[gray](FF) circle (0.35);
\fill[gray](GG) circle (0.35);

\fill[gray] (A)  circle (0.35);
\fill[gray] (4.7,-0.3) rectangle (5.3,0.3);
\fill[gray] (8.7,-0.3) rectangle (9.3,0.3);

\fill[black] (10.7,-0.3) rectangle (11.3,0.3);
\fill[gray] (12.2,-0.3) rectangle (12.8,0.3);

\fill[gray] (F)  circle (0.35);
\fill[gray] (G) circle (0.35);

\end{tikzpicture}}
\caption{Representation of a sample network in our framework.\label{fig:00}}
\end{center}
\end{figure}
\end{center}

\section{MINLP Formulation}\label{sec:2}

In this section we {provide} a Mathematical Programming formulation for the problem described in Section \ref{sec:1}. First of all, we {are given} an input network, like the one depicted in Figure \ref{fig:00}, including the topological route map and the stations, some of them being interchange stations and the lines pool $L$ defined over this network.  {In addition, short-turning decisions are allowed in some lines of the lines pool. These decisions always concern a set of given consecutive stations of those lines. In what follows and when no confusion is possible, we will refer to both, the set of consecutive stations and the trip concerned by a short-turning decision as a \textit{ short-turn}}. On the other hand, trips trough the full-length lines will be referred as \textit{whole trips}.  {In order to introduce the model we start by defining} the set of \emph{parameters} describing the remainder of the input data {as well as} the \emph{decision variables} used to identify a feasible solution.
\subsection{Parameters}

The input parameters for our model are the following:

\begin{itemize}
\setlength{\itemindent}{-0.2cm}
\item $[0,T]$: Planning horizon in which {trains start their journeys} at a given head of line station. Note that a train may start its last journey on $[0, T]$ while its last stop {may occur} in a time instant $t>T$.
\item $L$:  Set of lines in the network {(lines pool)}. As mentioned above, round-trip lines are considered as two different lines sharing the stations but {traversed} in opposite directions (rigorously speaking, the stations represent platforms of the corresponding line). Each one of the lines is assumed to be described by its node stations and its directed connections between consecutive stations. Additionally we will denote by $LS$ the set of lines containing a short-turn  and by $LNS$ the remainder {(those in which none of its proper subsets of stations can be activated as short-turns). Clearly, } $L=LS \cup LNS$.
\item $N_\ell=\{1, \ldots, n_\ell\}$: {Set of } stations of line $\ell\in L$. Stations are assumed to be ordered in its travelling direction. Observe that  if the lines $\ell$ and $\ell^\prime$ correspond to  {the same round-trip line but traversed} in opposite directions they have the same number of stations ($n_\ell=n_{\ell^\prime}$) and station $i\in\ell$ represents the opposite platform to station $n_{\ell^\prime}-i+1\in \ell^\prime$.

Additionally, for each line  $\ell\in LS$ we will denote by $S_\ell$ the set of stations in the (unique) short-turn,  being $1_{S_\ell}$ the first station of the short-turn and $n_{S_\ell}$ the last one.

In order to present a clearer Mathematical Programming model, we will also denote from now on by $\overline{\mathcal{S}}=\Big\{(i,\ell): i\in N_\ell, \ell\in LNS\Big\} \cup \Big\{(i,\ell): i\in N_\ell\backslash S_\ell, \ell\in LS\Big\}$, i.e., those stations which are not part of the available short-turns.
\item $d_i^\ell$: Distance{ (measured as travel time)}  between the stations $i$ and $i+1$ of line $\ell\in L$. {Recall, as mentioned in Section \ref{sec:1}, that} we assume that the speed of the trains operating between consecutive station $i$ and $i+1$ is fixed, { and thus, this distance is trip independent.}
\item ${e}_i^{\ell}$: Stopping time for  {any} train at station $i$ of line $\ell$ before leaving to station $i+1$. This value represents a time window to perform  different operations as the unload/load of passengers from/to the train or {the} fine-tuning of the train. Without loss of generality, we also consider that the stopping times are trip independent.
\item $t_{1\mapsto 1_{S_l}}$:   {Distance (measured as travel time plus stopping time at the last and intermediate stations) between the head of line station and the first station of the short-turn for}  line $l\in LS$. Note that, this parameter is given by the following expression:
     \begin{equation}
     t_{1\mapsto 1_{S_l}}=\sum_{r=1}^{1_{S_l}-1}({d_r^l}+e_{r+1}^{\ell}). \label{t1s}
     \end{equation}
\item $IS^\ell$: Minimum safety time interval between consecutive trips in a given line $\ell \in L$.
     \item $K_\ell= \{1, \ldots, {\bar{k}_\ell}\}$:  {Set of} trips in line $\ell\in L$.  It is worth noting that the exact number of trips in a line is a decision of the line planning process  and  which is not known beforehand. In fact,  some trips will not be actually used on the planning. We will refer to them as \textit{fake trips}. Note also that the maximum number of possible trips in the line $\ell\in L$ can always be upper bounded by  $\left\lceil\frac{T}{IS^\ell}\right\rceil+1$. Typically the value $\bar{k}_\ell$ will be smaller than this bound and should be estimated on the basis of the technical specifications of the network.
\item $Q = \{q_1, \ldots, q_{|Q|}\}$: Admissible capacities for trains operating in all the lines. In some cases, a single capacity profile is allowed while in some others a base capacity is available and, by adding extra wagons, it can be {doubled, triplicated}, etc.
\end{itemize}

We also need the proportion of passengers, which, not being able to catch the train in a given attempt, still insist on using the system and wait for the next train. In addition, our model corresponds to a line planning with route assignment in which the OD-matrix, whose entries describe the proportions of passengers moving among stations of each line is given.

\begin{itemize}
\setlength{\itemindent}{-0.2cm}
\item $\alpha$: Proportion of passengers that  cannot get on a train because of lack of capacity and await for the next one (proportion of persisting passengers). This parameter is assumed to be independent of stations and lines.  It can also model the probability  that a passenger gives up from getting on a crowded train and awaits for the next one.
\item $p^\ell_{ij}$: Proportion of passengers awaiting at station $i$ to go to the station $j$ using the line $\ell\in L$. Note that this parameter can assume positive values only if $i<j$, i.e., if station $i$ is previous to station $j$.
\end{itemize}

The following parameters will be used to assess a given planning solution. As explained in Section \ref{sec:1}, we propose the aggregation of the maintenance/operational planning costs together with a \emph{quantification} of the social cost whose purpose is to measure the quality of the service provided by the system. Obviously, the correct estimation of these costs is fundamental for a meaningful model.

\begin{itemize}
\setlength{\itemindent}{-0.2cm}
\item $b_q^\ell$: Fixed cost for starting a whole trip at line $\ell \in L$ with capacity $q\in Q$. Usually, the larger the capacities, the higher the fixed costs on the trips. These costs  model the consumption of resources to prepare a train for a new trip, the energy spent in the trip, the fixed costs of the staff needed to control the train, etc.
\item $b_{Sq}^\ell$: Fixed cost for starting a short-turn at line $\ell \in LS$ with capacity $q\in Q$. As in the previous cost, larger  capacities and lengths of the short-turns usually involve higher costs on the trips. These fixed costs are usually smaller than the corresponding ones for whole trips.
\item $\gamma^{{\ell}}_{ij}$: Unitary profit for transporting a passenger from the station $i$ to the station $j$ of the line $\ell\in L$. It represent the cost of the ticket paid by a single passenger to use the service for a given origin-destination trip.
\item $\mu_1$: Unitary penalty for passengers who  cannot get on the first arriving train due to its limited capacity and still insist on using the system. This parameter is an estimation of the loss supported by the service for unsatisfied passengers that still keep waiting to use the subway service.
\item $\mu_2$: Unitary penalty for passengers who give up using the network and leave the system after they cannot get on the first arriving train due to its limited capacity. This parameter is an estimation of the loss supported by the service for unsatisfied passengers that are lost by the system.
\end{itemize}

The above set of parameters are  listed in a more compact form in Table \ref{parameters}

\begin{table}[h]
\begin{center}
{\small\begin{tabular}{|cp{10cm}|}\hline
Parameter & Description\\\hline
$[0,T]$ & Planning horizon.\\
$L$ & Lines of the network {(Lines pool).}\\
$LS$ & Lines enabling short-turns.\\
$LNS$ & Lines without short-turns.\\
$N_\ell=\{1, \ldots, n_\ell\}$ & Stations of line $\ell\in L$. \\
$S_\ell= \{1_{S_\ell}, \ldots, n_{S_\ell}\}$ & Stations in the short-turn of line  $\ell\in LS$.\\
$\overline{\mathcal{S}}$ & Stations  {not affected by short-turns.} \\
$d_i^\ell$& {Distance (measured as time)} between stations $i$ and $i+1$ of line $\ell\in L$. \\
${e}_i^{\ell}$& Stopping time that a train spends in the station $i$ of the line $\ell\in L$.\\
$t_{1\mapsto 1_{S_l}}$ & {Distance (measured as time) between the head of line station and the first station of the short-turn  for} line $\ell\in LS$.\\
$IS^\ell$& Minimum safety time interval between consecutive trips, for  line $\ell\in L$.\\
$K_\ell= \{1, \ldots, {\bar{k}_\ell}\}$&  Set of trips in line $\ell\in L$, some of them being \textit{fake trips}. \\
$Q = \{q_1, \ldots, q_{|Q|}\}$& Admissible capacities for trains operating in all the lines.\\
$p^\ell_{ij}$ & Proportion of passengers moving between stations $i$ and $j$ of line $\ell\in L$. \\
$\alpha$ & Proportion of persisting passengers. \\
$b_q^\ell$ & Fixed cost for starting a whole trip at line $\ell\in L$ with capacity $q\in Q$.
\\$b_{Sq}^\ell$ & Fixed cost for starting a short-turn at line $\ell \in LS$ with capacity $q\in Q$.\\
$\gamma^{{\ell}}_{ij}$ & Unitary profit  {for transporting} a passenger  from the station $i$ to the station $j$ {of }line $\ell\in L$.\\
$\mu_1$ & Unitary penalty for persisting passengers.\\
$\mu_2$ & Unitary penalty for passengers who give up using the system.\\
\hline
\end{tabular}}
\end{center}
\caption{Parameters of the model.\label{parameters}}
\end{table}

\subsection{Decision variables}

The following set of decision variables are used in our Mathematical Programming model:

\begin{itemize}
\item {\bf Continuous Variables:}
\begin{itemize}
\setlength{\itemindent}{-0.2cm}
\item $t^{k{\ell}}_1$: Departure time from the initial station of line $\ell\in L $ at its $k$-th trip, $k\in K_\ell$.\\

Since the travel time between consecutive stations and the stopping time at each station are fixed, the departure time from the initial station of a line will be the reference time to calculate the departure time from the rest of stations of the line. 
\end{itemize}
It is worth mentioned  that  we will force any fake trip to operate on the line at the same time as the previous true trip. Thus, the departure time  from the initial station of the line of a fake trip must coincide with the departure time from the initial station of the previous true trip and therefore, the departure time from the rest of stations must also coincide.
Note also that  if $\ell\in LS $ some of its trips may be short-turns. When a short-turn is activated  the stations  not affected by the short-turn  are considered as nodes of a fake trip, being the departure times from these stations fixed to the ones of the previous true  whole trip. In order to avoid infeasible solutions due to inconsistent times relating departure times in consecutive stations (in and out the short-turn)  the following continuous variables are needed.
\begin{itemize}
\setlength{\itemindent}{-0.2cm}
\item  $w^{k{\ell}}$: Difference between the actual departure time from the first station of the short-turn of the $k$-th trip,  $k\in K_\ell$, of line $\ell\in LS$ and the time when it should depart from this station taking into account its departure time from the initial station of the line. Note that if  the $k$-th trip is a whole trip then, $w^{k{\ell}}=0$.\\
Note also that in the definition of these variables we are implicitly assuming that the initial station of the line is not part of the short-turn. If this occurs, we need to take as reference time the departure time from any other station out of the short-turn.
\end{itemize}

For modeling purpose  we need {to} distinguish between the flow of passengers that get on a train that performs a whole trip  and the flow of passengers that get on a train that only performs {a} short-turn. {Note that the first one is defined for any line of the lines pool whilst the second is only defined for those lines containing a short-turn}.
\begin{itemize}
\setlength{\itemindent}{-0.2cm}
\item $f^{k{\ell}}_i$: Flow of passengers captured in the station  {$i \in N_{\ell}$} by the train that covers the $k$-th trip of the line $\ell\in L$, {being $k\in K_\ell$  a whole trip}.
\item $g^{k{\ell}}_i$: Flow of passengers captured in the station $i\in S_\ell\backslash\{n_{S_\ell}\}$ by the train that covers the $k$-th trip of the line $\ell\in LS$, {being $k\in K_\ell$  a short-turn.}
\end{itemize}
Note that  the  $k$-th trip of a  line $\ell\in LS$  is either, a whole trip or a short-turn. Thus, if $g^{k\ell}_i>0$ then,  $f^{k\ell}_i=0$ and viceversa. Furthermore, if the $k$-th trip is a fake trip for line $\ell \in LS$ (resp. $\ell \in LNS$), then $g^{k\ell}_i=0$, for all $i\in S_\ell\backslash\{n_{S_\ell}\}$ and $f^{k\ell}_i=0$ for all $i\in N_\ell$ (resp. $f^{k\ell}_i=0$ for all $i\in N_\ell$).

\item {\bf Binary Variables:}
\begin{itemize}
\setlength{\itemindent}{-0.2cm}
\item The decisions about the capacities of the trains at each trip are modeled using the following variables:
$$
y_q^{k\ell}  = \left\{\begin{array}{cl}1 & \mbox{if the $k$-th trip of line $\ell$ is {a whole trip}}\\  &\mbox{with capacity $q$ }\\
0 & \mbox{otherwise,}\end{array}\right.
\quad {k \in K_\ell},\, q\in Q,\, \ell\in L.$$
$$
y_{Sq}^{k\ell}  = \left\{\begin{array}{cl}1 & \mbox{if the $k$-th trip of line $\ell$ traverses every } \\
& \mbox{station of the short-turn  $S_\ell$ with capacity $q$}\\
0 & \mbox{otherwise,}\end{array}\right.
\quad {k \in K_\ell},\, q\in Q,\, \ell\in LS.$$
\end{itemize}
Note that  {if the $k$-th trip of  line $\ell \in LS$ is a short-turn with capacity $q$ then, } $y_{Sq}^{k\ell}=1$ and  $y_q^{k\ell}=0$. On  {the other hand}, {if the $k$-th trip of  line $\ell \in LS$ is a whole trip with capacity $q$ ($y_q^{k\ell}=1$)}, the train traverses all the stations of the line and, in particular, the stations of the short-turn, and  then, $y_{Sq}^{k\ell}=1$. So, we have the following relationship between the two sets of variables:
 \begin{equation}\label{rely}
y_q^{k\ell}\leq y_{Sq}^{k\ell}\quad {k \in K_\ell},\, q\in Q,\, \ell\in LS.
\end{equation}
{Observe that we can detect whether the $k$-th trip is a fake trip} {for line $\ell \in LS$ (resp. $\ell \in LNS$)} {by checking if $\sum_{q \in Q} y_{Sq}^{k\ell} =0 ( =\sum_{q \in Q} y_{q}^{k\ell})$}
{(resp. $\sum_{q \in Q} y_{q}^{k\ell}=0$).}

\item {\bf Auxiliary Variables}: A set of auxiliary variables, computed from the previous ones, is considered in order to ease the reading of our Mathematical Programming formulation.
\begin{itemize}
\setlength{\itemindent}{-0.2cm}
\item $t^{k\ell}_i$: Time instant in which a train departs from station $i$ at the $k$-th trip of line $\ell\in L$. This value can be obtained by adding the travel times  between consecutive stations plus the stopping times at the traversed stations. Thus, this time is given by the following expressions:
 \begin{align}
&t^{k\ell}_i  =t_1^{k\ell}+\sum_{r=1}^{i-1}({d_r^\ell}+e_{r+1}^{\ell}), \quad i>1,\, (i,{\ell})\in \overline{\mathcal{S}},\,  k \in K_\ell,\label{tiemposalida1}\\
& t^{k\ell}_i=t_1^{k\ell}+\sum_{r=1}^{i-1}({d_r^\ell}+e_{r+1}^{\ell})+w^{k{\ell}},  \quad i>1,\, i\in  S_\ell,\,  k \in K_\ell,\, \ell \in LS.\label{tiemposalida3}
 \end{align}
Expressions \eqref{tiemposalida1} allow us to appropriately define our auxiliary variables $t_i^{k\ell}$ by updating the time in which a train starts its journey {at the first station of the line} to the time spent until it leaves station $i$ for stations outside a short-turn. Equations \eqref{tiemposalida3}  allow us to compute the  arrival time of a train to any  station of the short-turn.  As mentioned before, trips that only cover stations in the short-turn  are considered as fake trips for the stations which are not part of the short-turn and their corresponding departure times are supposed to be the same as the times of the previous trip.  Variables  $w^{k{\ell}}$ permit to fit departure times of the stations in each short-turn trip. These variables will {take value} $0$ if the train cover the complete line.
  \item  $D_i^\ell(t^{k\ell}_i)$: Number of passengers accumulated at  station $i$  of line $\ell\in L$ {from the origin of the planning horizon} up to  time $t^{k\ell}_i$.  This number is given by a function  of the time  $t^{k\ell}_i$  and takes into account the external arrivals and also the passengers arriving to the station to change to another line. As mentioned {above}, a complete description of this \textit{demand function} is given in Section \ref{sec:3}.

\item $h^{k\ell}_i$: Excess of passengers that where not able to get on the train at station $i$ of the $k$-th trip of line $\ell\in L$ because a lack of capacity of the train. To compute this variable we distinguish between whole trips  in any  of the pool of lines and short-turn trips in those lines containing a short-turn. In the last case we assume that the passengers catching a short-turn train with destination to a station outside the short-turn get-off from the train in the last station of the short-turn to catch a {whole}  trip train of the same line. Taking into account these assumptions, the excess of passengers can be computed as follows:
          \begin{align}
& {h^{1\ell}_i} = D_i^\ell(t_i^{1\ell}) - f_i^{1\ell},\quad \text{for }  (i,{\ell})\in \overline{\mathcal{S}},\label{remanente0a}\\
& {h^{1\ell}_i} = D_i^\ell(t_i^{1\ell}) - f_i^{1\ell} - g_i^{1\ell}, \quad \text{for }  i\in S_\ell\backslash\{n_{S_\ell}\},\, \ell \in LS,\label{remanente0c}\\
& h^{1\ell}_i = D_i^\ell(t_i^{1\ell}) - f_i^{1\ell} +  \dsum_{r=1_{S_\ell}}^{i-1} \dsum_{j=i+1}^{n_\ell} p_{rj} g_r^{1\ell},\quad \text{for }  i=n_{S_\ell},\, \ell \in LS,\label{remanente0d}\\
& h^{k\ell}_i =D^\ell_i(t^{k\ell}_i)-D^\ell_i(t^{(k-1)\ell}_i)+\alpha h^{(k-1)\ell}_i-f_{i}^{k\ell}, \quad \text{for }  (i,{\ell})\in \overline{\mathcal{S}}, \, k=2,\ldots, \bar{k}_\ell, \label{remanentea}\\
& h^{k\ell}_{i} = D^\ell_i(t^{k\ell}_i)-D^{{\ell}}_i(t^{(k-1)\ell}_i)+\alpha h^{(k-1)\ell}_{i} - f_{i}^{k\ell}- g_i^{k\ell}, \label{remanenteS}\nonumber\\
& \hspace*{5cm} \text{for } i \in S_\ell \backslash\ \{n_{S_\ell}\}, \;  k=2,\ldots, \bar{k}_\ell, \ell \in LS, \\
& h^{k\ell}_i =  D^\ell_i(t^{k\ell}_i)-D^{{\ell}}_i(t^{(k-1)\ell}_i)+\alpha h^{(k-1)\ell}_{i} - f_{i}^{k\ell} +  \dsum_{r=1_{S_\ell}}^{i-1} \dsum_{j=i+1}^{n_\ell} p_{rj} g_r^{k\ell},\label{remanented} \nonumber\\
& \hspace*{6cm} \text{for }  i=n_{S_\ell}, k=2,\ldots, \bar{k}_\ell, \ell \in LS.
      \end{align}
In the first trip, the excess of passenger is {computed} by subtracting to the demand the flow of passengers already caught by the train ({equations} \eqref{remanente0a}-\eqref{remanente0d}). For modeling the flow of caught passengers we take into account in \eqref{remanente0c} that the trip could be either, a whole trip  ($f_i^{1\ell}$) or a short-turn ($g_i^{1\ell}$). In the special case of the last station of the short-turn, one needs to add to the excess of passengers, those that get-off from a {train of a}  short-turn {trip} to catch a train of a whole trip of the same line (equations \eqref{remanente0d}).
For the rest of the trips ({equations} \eqref{remanentea}--\eqref{remanented}), we take into account the demand accumulated after the previous trip plus the excess of persisting passengers of the previous trip.
\end{itemize}

\item {\bf Semicontinuous  Variables}:
  We also consider a {set of} semicontinuous variable{s} collecting the excess of passengers only for \textit{true trips}:
$$
x^{k\ell}_i=\left\{\begin{array}{ll}
 h_{i}^{k\ell} & \text{ if $k$ is a true trip for station $i$ of line $\ell$}\\
 0 &\text{otherwise,}
 \end{array}\right. \quad {k \in K_\ell},\, i\in N_\ell,\, \ell\in L.
 $$
This set of variables will allow us to account in the objective function for the actual excess of passengers, i.e., the ones associated to true trips.
\end{itemize}

Table \ref{variables} summarizes the set of {above mentioned} variables.

\begin{table}[h]
\begin{center}
{\small\begin{tabular}{|cp{13cm}|}\hline
Variable & Description\\\hline
$t^{k{\ell}}_1$ & Departure time from the initial station of line $\ell\in L$ at its $k$-th trip.\\
$f^{k{\ell}}_i$ & Flow of passengers captured in the station $i$ by the train that covers the $k$-th trip of the line $\ell\in L$, when $k$ {is a whole trip.}\\
$g^{k{\ell}}_i$ & Flow of passengers captured in the station $i\in S_\ell \backslash{\{n_\ell\}}$ by the train that covers the $k$-th trip of the line $\ell\in LS$, when $k$ only covers the short-turn stations.\\
$w^{k{\ell}}$ &  {Difference between the actual departure time from the first station of the short-turn of the $k$-th trip of line $\ell\in LS$ and the time when it should depart from this station taking into account its departure time from the initial station of the line.}\\
$y_q^{k\ell}$   & $\left\{\begin{array}{cl}1 & \mbox{if the $k$-th trip of line $\ell\in L$ {is a whole trip}  with capacity $q$}\\
0 & \mbox{otherwise}\end{array}\right.$\\
$y_{Sq}^{k\ell}$ & $\left\{\begin{array}{cl}1 & \mbox{if the $k$-th trip of line $\ell\in LS$ {traverses the short-turn  $S_\ell$} with capacity $q$.}\\
0 & \mbox{otherwise}\end{array}\right.$\\
{$t^{k\ell}_i$} & {Departure time from the station $i$ of line $\ell$ in its $k$-th trip.}\\
$D_i^\ell(t^{k\ell}_i)$ & Number of passengers accumulated from instant $0$ up to instant $t^{k\ell}_i$ in the station $i$ of line $\ell\in L$.\\
$h^{k\ell}_i$& Excess of passengers that where not able to get on the train at station $i$ at the $k$-th trip of line $\ell\in L$ because of a lack of capacity.\\
 $x^{k\ell}_i$ &  Excess of passengers only if $k$ is a \textit{true} trip for station $i$ of line $\ell\in L$.\\
 \hline
\end{tabular}}
\end{center}
\caption{Variables of the model.\label{variables}}
\end{table}

Using the variables described above, we give now the formulation of our problem  distinguishing the two  main elements: (1) an objective function aggregating the economical and social costs of any feasible solution; and (2) a set of technical constraints where the line concept and the timetable are specified in order to serve the estimated flow of passengers over the planning horizon.

\subsection{Objective Function}

In what follows we describe the different {\bf costs} and {\bf rewards} that are considered in the objective function of our model, for a given line $\ell{\in L}$, in terms of the variables and parameters described above.

\begin{itemize}
\setlength{\itemindent}{-0.2cm}
\item {\bf Capacity cost: }\\
It accounts for the costs $b_q^\ell$ and $b^\ell_{Sq}$ corresponding to the train capacity $q\in Q$ for every trip of  line $\ell \in L$. This capacity is controlled by means of the variables $y^{k\ell}_q$  or $y^{k\ell}_{Sq}$ depending {on} whether the line contains short-turns or not.
\begin{equation}
\left\{
\begin{array}{cl}
\dsum_{k \in K_\ell} \dsum_{q\in Q}  b^\ell_q y^{k\ell}_q &\mbox{ if } { \ell\in LNS},\\
{\dsum_{k \in K_\ell} \dsum_{q\in Q}  b^\ell_q y^{k\ell}_q  + }\dsum_{k \in K_\ell} \dsum_{q\in Q}  b^\ell_{Sq} (y^{k\ell}_{Sq} - y^{k\ell}_{q}) & \mbox{ if } \ell\in LS.
\end{array}\right.
\label{obj:2}\tag{{${\rm Cap}(\ell)$}}
\end{equation}

{Observe that in the second expression, if proper short-turns {with capacity $q$} are activated, one has $y_{Sq}^{k\ell}=1$  and $y_{q}^{k\ell}=0$,  {and thus,} the amount $b_{Sq}^\ell$ is accounted {for}. In case of whole trips {with capacity $q$}, by \eqref{rely}, one has that $y_{Sq}^{k\ell}=1$ and $y_{q}^{k\ell}=1$ . Thus, the second addend is zero, and only $b_{q}^\ell$ is accounted {for}, modeling adequately the capacity costs.}

\item {\bf Reward per served passengers:} \\
{The unitary reward per served passenger} computes the estimated revenue received when  passengers use a line. It  is obtained as {the} average revenue given by the different transport tickets used by the passengers. In order to compute {the overall reward per served passengers}, observe that {if the $k$-th trip is a whole trip,} the expression $ \sum_{r=0}^{i-1} p_{ri} f_r^{k\ell}$   returns the estimated number of passengers getting off at the station $i$ (coming from any other previous station). Note that, in case a short-turn {is} activated on {any  line $\ell\in LS$}, it is needed to add the rewards of passengers being routed within the short-turn {($\sum_{i \in S_\ell}  \sum_{r=1_{S_\ell}}^{i-1}  p_{ri}^\ell g_r^{k\ell} $)} together with the reward from those which use the short-turn to get-off at the last station and continue {with the next whole trip}.  Then, the overall reward per served passengers is:
\begin{equation}
\hspace*{-1cm}
\left\{\begin{array}{cl}
\dsum_{i \in N_\ell} \dsum_{k\in K_\ell} \dsum_{r=0}^{i-1}\gamma^{\ell}_{ri} p_{ri}^\ell f_r^{k\ell} & \mbox{ if } { \ell\in LNS},\\
\dsum_{k\in K_\ell}\left(\dsum_{i \in N_\ell}  \dsum_{r=0}^{i-1}\gamma^{\ell}_{ri} p_{ri}^\ell f_r^{k\ell}+\dsum_{i \in S_\ell}  \dsum_{r=1_{S_\ell}}^{i-1}\gamma^{\ell}_{ri}   p_{ri}^\ell g_r^{k\ell} + \hspace*{-0.2cm}\dsum_{r \in S_\ell:\atop
 r\neq n_{S_\ell}} \dsum_{j=n_{S_\ell+1}}^{n_\ell}\gamma^{\ell}_{rn_{S_\ell}}   p_{rj}^\ell g_r^{k\ell}\right)& \mbox{ if } \ell\in LS.
\end{array}\right.
\label{obj:4a}\tag{{${\rm RewPPass}(\ell)$}}
\end{equation}

\item {\bf Costs of non-served passengers:} \\
It accounts for the social cost incurred when a passenger cannot get on the train arriving to the station due to its lack of capacity. The unitary cost should aggregate some indicators of the service quality and some subjective measures of the satisfaction degree perceived by the passengers. The overall cost is computed by using the total number of passengers exceeding the capacity of the system at some instant of the planning horizon, scaled by the average penalties of persisting/giving up passengers, as follows:
\begin{equation}\alpha \mu_1\sum_{i\in N_\ell}
\sum_{k\in K_\ell} x_i^{k\ell}+  (1-\alpha) \mu_2 \sum_{i\in N_\ell} \sum_{k\in K_\ell} x_i^{k\ell}, \label{obj:5}\tag{{${\rm NonServed}(\ell)$}}
\end{equation}
Observe that a high excess of passenger at the end of the planning horizon can be easily avoided by increasing the $\mu$-parameters for the last trip or by adding appropriate constraints involving the $x$-variables.
\end{itemize}

The overall cost of using line $\ell\in L$ during the planning horizon can be expressed as:

\begin{equation}
\eqref{obj:2} -\eqref{obj:4a} +\eqref{obj:5}\label{obj}\tag{{\rm COST}($\ell$)}
\end{equation}

\subsection{Modelling Constraints}

{In what follows we describe the constraints linking the variables and parameters in our model. They have been classified in four main blocks: capacity constraints, time control constraints, flow control constraints and passenger surplus constraints.}

\begin{itemize}

 \item Capacities and true/fake trips:
          \begin{subequations}
    \makeatletter
        \def\@currentlabel{${\rm C}1$}
        \makeatother
       \label{CS1}
        \renewcommand{\theequation}{${\rm C}1-{\arabic{equation}}$}
       \begin{align}
    \sum_{q\in Q} y^{1\ell}_q=  1, &\qquad\ell\in LNS,\label{c:5}\\
    \sum_{q\in Q} y^{k\ell}_q\leq 1, &\qquad k>1, \ell\in L,\label{c:6}\\
     \sum_{q\in Q} y^{\overline{k}_\ell \ell}_q=  1, &\qquad\ell\in L \label{lastround},\\
     y_q^{k\ell}    \leq y_{Sq}^{k\ell} ,& \qquad q\in Q,  k  \in K_\ell, \ell\in LS,\label{c:ys}\\
\sum_{q\in Q} y^{1\ell}_q+\sum_{q\in Q} y^{1\ell}_{Sq}\geq  1,&\qquad  \ell\in LS,\label{c:5b}\\
 \sum_{q\in Q} y^{\kappa_\ell \ell}_q=\sum_{q\in Q} y^{1\ell}_{Sq}- \sum_{q\in Q}y^{1\ell}_q,&\qquad\ell\in LS,\label{sub2}\\
\sum_{q\in Q} y^{k\ell}_{Sq} \leq 1,&\qquad k>1, \ell\in LS,\label{c:6b}
   \end{align}
    \end{subequations}
where $\kappa_\ell=\left[\frac{t_{1\mapsto 1_{S_l}}}{IS^\ell}\right]+1$, i.e. the number of short-turn trips of line $\ell {\in LS}$ which fit within the period of time taken by a train to go from the head of line to the first station of the short-turn.

When short-turns are not allowed for a line, the appropriate definition of the  capacity variables is ensured by constraints \eqref{c:5}--\eqref{lastround}. They enforce that exactly one {of the allowed capacities}  is chosen for the first and the last trip and at most one for the rest of them. Fake (resp. true) trips are identified by trips with capacities equal to (resp. greater than) zero. Thus,  constraints \eqref{c:5} and \eqref{lastround}  determine that the first and the last trip of each line are true trips. We will see later that this permits the actual trains to be scheduled from the beginning to the end of the planning horizon, providing the users a complete service during that time interval. Note that  constraints \eqref{c:6} and \eqref{lastround} are also valid for lines {allowing} short-turns  and then, when {short-turns are allowed for a line}, the appropriate definition of the capacity variables is warranted by constraints \eqref{c:6}--\eqref{c:6b}.   Constraints \eqref{c:ys} indicate that {when a whole trip is a true trip, it is also a true trip for the  short-turn stations}. Constraints \eqref{c:5b} fix that the first trip ({being either, a whole trip or a short-turn}) is a true trip. Constraints \eqref{sub2} force trip $\kappa_\ell$ to be a true {whole} trip  (resp. a fake  trip) if the first trip {is a }  short-turn (resp. if the first trip {is a whole trip}).
These constraints, together with \eqref{c:5b} are the equivalent to \eqref{c:5} for lines with short-turns, and they ensure that a real train is scheduled from the beginning of the planning horizon. Finally, constraints \eqref{c:6b}  are the equivalent to \eqref{c:6} for short-turn trips.

\item Time control:

       \begin{subequations}
    \makeatletter
        \def\@currentlabel{${\rm C}2$}
        \makeatother
       \label{CS2}
        \renewcommand{\theequation}{${\rm C}2-{\arabic{equation}}$}
       \begin{align}
 t_1^{1\ell}=0,&\qquad  \ell\in L,\label{firsttime}\\
     t_1^{\overline{k}\ell}=T, &\qquad \ell\in L,\label{lasttime}\\
    t_1^{\kappa_\ell \ell}\leq T\left(1-\sum_{q\in Q} y^{1\ell}_{Sq}+\sum_{q\in Q} y^{1\ell}_q \right),&\qquad \ell\in LS,\label{sub1}\\
    IS \left(\dsum_{q \in Q} y_q^{k\ell}\right) \leq t^{k\ell}_i- t^{(k-1)\ell}_i, &  \label{c:8a} \nonumber\\
     \hspace*{4cm} k=2,\ldots, \overline{k}_\ell, &\, {(i,l)\in \overline{\mathcal{S}} \mbox{ with $k\neq \kappa_\ell$ if  $\ell\in LS$,}}\\
 t^{k\ell}_i- t^{(k-1)\ell}_i \leq T\left(\sum_{q\in Q} y^{k\ell}_q\right),& \quad    k=2,\ldots, \overline{k}_\ell,\, {(i,l)\in \overline{\mathcal{S}}} \label{c:8b}\\
 IS \left(\dsum_{q \in Q} y_{Sq}^{k\ell}\right) \leq t^{k\ell}_i- t^{(k-1)\ell}_i, &\quad   i \in S_\ell, \, k=2,\ldots, \overline{k}_\ell, \, \ell\in LS,\label{c:8c}\\
 t^{k\ell}_i- t^{(k-1)\ell}_i \leq (T+t_{1\mapsto 1_{S_l}})\left(\sum_{q\in Q} y^{k\ell}_{Sq}\right),&
 \quad i \in S_\ell,k=2,\ldots, \overline{k}_\ell, \, \ell\in LS, \label{c:8d}\\
    - t_{1\mapsto 1_{S_l}}  \left(1-\sum_{q\in Q_\ell} y_{q}^{k\ell}\right) \leq  w^{k\ell}, & \quad
k \in K_\ell,\, \ell \in LS,\label{tiemposalida:w1}\\
    w^{k\ell}\leq (T+t_{1\mapsto 1_{S_l}}) \left(1-\sum_{q\in Q_\ell}  y_{q}^{k\ell}\right), &\quad
k \in K_\ell,\, \ell \in LS, \label{tiemposalida:w2}
\end{align}
    \end{subequations}

Constraints \eqref{firsttime} and \eqref{lasttime} state that the first and the last trip of each line should exactly start at the first station of the complete line at instant time $0$ and  $T$, respectively. If the line does not contain a short-turn, these constraints together with  \eqref{c:5} and \eqref{lastround}  ensure that there are trains traveling the line during the whole planning horizon. If the line contains short-turns, recall that constraints \eqref{c:5b} permit the first trip to be a short-turn trip. In this case, constraints \eqref{sub1} together with constraints \eqref{sub2} force trip $\kappa_\ell$ to be a true {whole} trip  starting at time $0$ at the head of the line.

Constraints \eqref{c:8a} ensure that the arrival times between consecutive trains satisfy the safety time window interval.  Constraints \eqref{c:8b} force a fake trip to operate on the line at the same time as the previous true trip.  For lines allowing short-turns constraints \eqref{c:8a}--\eqref{c:8d}  represent the same as the above but taking into account that  in this case the trip can be either, a whole trip or a short-turn.

Observe that in \eqref{c:8a}  if $\ell\in LS$ the case $k=\kappa_\ell$ is excluded. The reason is that constraints  \eqref{sub2}  force trip $\kappa_\ell$ to be a true whole trip  if and only if the first trip ($k=1$) is a short-turn and, in this case, constraints \eqref{sub1} ensure that trip $\kappa_\ell$ starts at the head of the line station at time 0, and then, it is the first whole trip. Finally, constraints \eqref{tiemposalida:w1} and \eqref{tiemposalida:w2} fix the upper and lower bounds on the values of variables $w^{kl}$  when the $k$-th trip is a short-turn, enforcing that this variable is $0$ if the  trip is a whole trip.

     \item Flow control:

             \begin{subequations}
    \makeatletter
        \def\@currentlabel{${\rm C}3$}
        \makeatother
       \label{CS3}
        \renewcommand{\theequation}{${\rm C}3-{\arabic{equation}}$}
    \begin{align}
f_{i}^{k\ell} + \dsum_{r=1} ^{i-1} f_r^{k\ell}  \left(\dsum_{j=i+1}^{n_\ell} p_{rj}^\ell\right) \leq \dsum_{q \in Q} q y_q^{k\ell},& \quad k\in K_\ell,  i\in N_\ell,  \ell \in L ,\label{c:14d}\\
 g_{i}^{k\ell} + \dsum_{r=1_{S_\ell}} ^{i-1} g_r^{k\ell} \left(\dsum_{j=i+1}^{n_{S_\ell}} p_{rj}^\ell\right) \leq \dsum_{q \in Q} q (y_{Sq}^{k\ell}- y_{q}^{k\ell}),&\;  k\in K_\ell, i \in S_\ell\backslash\{n_{S_\ell}\},  \ell \in LS,\label{c:14i}\\
f_{i}^{1\ell}  {\leq} D_i^\ell(t_i^{1\ell}), &\quad i \in N_\ell, \ell \in L,\label{c:14b1}\\
f_{i}^{kl}\leq D^\ell_i(t^{k\ell}_i)-D^\ell_i(t^{(k-1)\ell}_i)+\alpha h^{(k-1)\ell}_i, &\quad  k=2,\ldots, \overline{k}_\ell, i\in N_\ell, \ell \in L,\label{c:14a1}\\
 g_{i}^{1\ell}\leq  D^\ell_i(t^{1\ell}_i), &\;\;i \in S_\ell\backslash\{n_{S_\ell}\}, \ell \in LS,\label{c:14b2}\\
  g_{i}^{k\ell}\leq  \left(D^\ell_i(t^{k\ell}_i)-D^\ell_i(t^{(k-1)\ell}_i)\right)+\alpha h^{(k-1)\ell}_{i},& \;\; k=2,\ldots, \overline{k}_\ell, i \in S_\ell\backslash\{n_{S_\ell}\},  \ell \in LS.\label{c:14a2}
 \end{align}
    \end{subequations}

The flow of passengers catching a given train is determined by the capacity of the train and the {mobility pattern} of people. Hence, the effective capacity of the trains arriving to a given station depends of the passengers that caught the train in previous stations of this same trip and whose destination is a subsequent station of the line. Such an effective capacity is warranted, depending on the case (short-turn allowed or not),  by constraints  \eqref{c:14d} and \eqref{c:14i}.
With constraints \eqref{c:14b1} (resp. \eqref{c:14a1}) we ensure that the flow of passengers captured at a given station by a train that covers the first trip (resp. the $k$-th trip for $k>1$) is at most the demand of passengers accumulated at station $i$ since the beginning of the planning horizon (resp. since the instant in which the previous train departed from that station plus the passengers that were not able to get on the previous train because of lack of capacity and wait for the next one). Constraints \eqref{c:14b2} and \eqref{c:14a2} are the {analogous} ones for short-turn {trips}.

\item Passenger surplus:

In order to compute only the surplus of passengers of a true trip, we use the  {set of semi-continuous variables $x^{k\ell}_i$}:
$$
x^{k\ell}_{i} = \left\{\begin{array}{cl}
h^{k\ell}_{i}  \times \sum_{q\in Q} y_{q}^{k\ell} & \mbox{ if $(i,l)\in \overline{\mathcal{S}}$ { or  ($i=n_{S_\ell}$, $\ell \in LS$)}},\\
h^{k\ell}_{i}  \times \sum_{q\in Q} y_{Sq}^{k\ell} & \mbox{ if $i\in S_\ell\backslash\{n_{S_\ell}\}$, $\ell \in LS$,}
\end{array}\right.
$$
which can be linearized as follows:
                      \begin{subequations}
    \makeatletter
        \def\@currentlabel{${\rm C}4$}
        \makeatother
       \label{CS4}
        \renewcommand{\theequation}{${\rm C}4-{\arabic{equation}}$}
            \begin{align}
& x^{k\ell}_{i}\geq  h_{i}^{k\ell}-M_i^\ell\left(1-\sum_{q\in Q_\ell} y_{q}^{k\ell}\right),  (i,l)\in \overline{\mathcal{S}} \mbox{ {or ($i=n_{S_\ell}$, $\ell \in LS$)}} ,\label{remanente0ax1}\\
& x^{k\ell}_{i}\geq  h_{i}^{k\ell}-M_i^\ell\left(1-\sum_{q\in Q_\ell} y_{Sq}^{k\ell}\right),  i\in S_\ell\backslash\{n_{S_\ell}\}, \ell \in LS,\label{remanente0ax3}
      \end{align}
    \end{subequations}
being $M_i^\ell$ a large enough constant bounding the surplus of passengers at any station $i$ of line $\ell\in L$.

\end{itemize}

\subsection{A compact MINLP formulation}

According to the decision variables{, the objective function} and the constraints described {above},  {the following Mathematical Programming formulation is valid for our {line} planning and timetabling model:}
\begin{align}
 \qquad\min & \dsum_{\ell\in L}\mbox{\ref{obj}} & \nonumber\\
\mbox{s.t. }  & \eqref{CS1},\, \eqref{CS2},\, \eqref{CS3} \mbox{ and } \eqref{CS4},\nonumber\\
  & 0  \le t^{k\ell}_1 \leq T,&\quad k\in K_\ell,\, \ell\in L, \nonumber\\
&f_i^{k\ell} \geq 0,  &\quad  k\in K_\ell,\, i\in N_\ell,\, \ell\in L,  \label{MINLP}\tag{P}\\
&g_i^{k\ell}\geq  0, &\quad  k\in K_\ell,\, i\in S_\ell\backslash\{n_{S_\ell}\},\,  \ell\in LS, \nonumber\\
&w^{k\ell} \in \mathbb{R}, &\quad k\in K_\ell,\, \ell \in LS, \nonumber\\
&x^{k\ell}_{i} \geq 0, &\quad k\in K_\ell,\, i\in N_\ell,\, \ell \in L,\nonumber\\
&y^{k\ell}_q  \in \{0,1\}, &\quad  k\in K_\ell,\, q\in Q,\, \ell\in L, \nonumber\\
& y^{k\ell}_{S_q}  \in \{0,1\}, &\quad  k\in K_\ell,\, q\in Q, \ell \in LS.\nonumber
\end{align}
Observe that although the above formulation seems to be separable by lines in $L$, the lines are linked through the demand function{{s}} $D_i^\ell(t)$ (constraints \eqref{c:14b1}-\eqref{c:14a2}) which represent the accumulated flow of passengers awaiting for a train at a given station $i$ of line $l\in L$ at time instant $t$. As we will describe in Section \ref{sec:3}, such a flow is affected not only by the line $\ell$ but also by other lines through passengers changing of lines at transfer stations. This function introduces new variables and non linear constraints into the above formulation.

{Several extensions may be easily accommodated within the above model as highlighted in the following remarks:}
\begin{remark} \label{remark1}
In our model the speed of trains is considered to be constant during the whole journey. However, one can easily modify expressions {\eqref{t1s}, }\eqref{tiemposalida1} and \eqref{tiemposalida3} using variables $v_{i}^{k\ell}{>0}$ to decide the speed of the train during a trip $k$ of line $\ell$ between stations $i$ and $i+1$. For instance, let $\rho_i$ be the physical distance between stations $i$ and $i+1$ and let $\omega_{i}^{k\ell}$ represent the inverse of the speed, i.e., $\omega_{i}^{k\ell} = \frac{1}{v_i^{k\ell}}$, then one could replace the travel times $d_i^\ell$ by $\rho_i \times \omega_{i}^{k\ell}$ in expressions {\eqref{t1s}, } \eqref{tiemposalida1} and \eqref{tiemposalida3}. By adding to the objective function a cost assessing the resource consumption due to speed changes one can have a more general model preserving the structure of the one stated above. Similar modifications can be also considered by enabling the model to decide about variable stopping times at any station.
\end{remark}

\begin{remark}\label{remark2}
The model allows us to join two consecutive planning horizons by passing data about numbers of passengers and arrival/departure times obtained from an optimal solution on the first planning period as input data for the {second} one{, and so on}. In particular, the passengers that may remain at station $i$ of  line $\ell\in L$ at the end of the first planning horizon can be considered as passengers at station $i$ to use line $\ell$ at the beginning of the {second} planning horizon{, and so on}. This information has to be incorporated in order to compute the demand function, as we will see in the next section.
\end{remark}

\section{The Demand function}
\label{sec:3}

One of the main goals of our model is to incorporate, in the design of the line planning and timetabling of an existing network, information about the flow of passengers moving through the network during the planning horizon. Clearly, the flow of passengers arriving to a given station is a random variable. Thus, we will incorporate to the model an estimation of its average value.

In order to model the number of passengers entering to the transportation system through a given station of a fixed line, we use the so-called \textit{demand function}, which maps at a given instant $t$ the accumulated number of passengers wanting to catch a train at this station (from the beginning of the planning horizon). Here, the estimation process should be carefully done in order to capture the essential behaviour of the demands served by the system.

{Different shapes for the function are possible within this framework to approximate the demand. The choice of such a shape is a crucial step in the modeling process since one has to find an equilibrium between obtaining accurate estimations and providing manageable mathematical programming formulations. Once again, motivated by our pilot experience with Metrolab{\textregistered}, we use a piecewise linear approximation whose slope is fixed for any given station $i \in N_\ell$, but whose breakpoints and discontinuities may change according to the flow induced by  external block of arrivals and by the rest of the lines. For a given line $\ell \in L$ and a station $i\in N_\ell$, we estimate the demand function as follows:
\begin{equation}\label{demand}\tag{${\rm D}$}
D_i^\ell(t) = \beta_{0i}^\ell +  \beta_i^\ell t + J^E_{i\ell} (t) + \dsum_{\ell'\neq \ell, \ell' \ni i} J^I_{i\ell\ell'} (t),
\end{equation}
for $t \in [0,\widehat{T}_\ell]$, with $\widehat{T}_\ell=T+\sum_{r=1}^{n_\ell-1} (d_r^\ell+e_{r+1}^\ell)$, {i.e., the maximum time in which the train can reach the last station of the line,}  and where:
\begin{itemize}
\setlength{\itemindent}{-0.8cm}
\item[] $\beta_{0i}^\ell$ is the number of passengers awaiting a train of line $\ell\in L$ in the station $i\in N_\ell$ at the beginning of the planning horizon.
\item[] $\beta_i^\ell$ is the average rate of passengers arriving to the station $i\in N_\ell$ of line $\ell\in L$ by unit of time.
\item[] $J^E_{i\ell}(t)$ is the sum of the external block of arrivals of passengers up to the instant $t$ to the station $i\in N_\ell$ of line $\ell\in L$.
\item[] $J^I_{i\ell\ell'} (t)$ is the sum of the block arrivals of passengers up to the instant $t$ to the station $i\in N_\ell$ of line $\ell\in L$ from line $\ell^\prime\in L$.
\end{itemize}
In what follows we will refer to $[0,\widehat{T}_\ell]$  as the \textit{extended planning horizon for line $\ell$}.
Thus, the demand function at a given time instant $t$ in station $i\in N_\ell$ of line $\ell\in L$, consists of three parts. The first one is a linear part, in which, from an initial number of passengers, $\beta_{0i}^\ell$, the number increases by a rate $\beta_i^\ell$.  However, such a base estimation may be modified either by external block of arrivals  (second part), $J^E_{i\ell}(t)$, or by passengers coming from other interacting lines at interchange stations (third part), $J^I_{i\ell\ell'} (t)$. 

As can be seen, in the formulation \eqref{MINLP}, the demand function $D_i^\ell(t)$ is used exclusively to access flows in the set of instants $t = t_{i}^{k\ell}$ for $i \in N_l, \ell \in L$ and $k \in K_\ell$. Each of the time instants in which the demand function needs to be evaluated induces {some} sets of inequalities and variables as those described in subsections \ref{ss:EA} and \ref{ss:IA} (for external and internal arrivals).

\subsection{External Arrivals: $J^E$\label{ss:EA}}

We consider that we are given both a set of breakpoints representing time instants when the block of arrivals occur and the amounts of passengers entering to the system at these instants for each station $i \in N_\ell$. {That is, we assume that a set of sorted instants $se^{i\ell}_1 < \cdots < {se_{re^{i\ell}}^{i\ell}}$ as well as discontinuity flow jumps associated to each of those instants $\Psi^{i\ell}_{1}, \ldots,  {\Psi^{i\ell}_{re^{i\ell}}}$ are known}, {i.e., a block arrival of  $\Psi^{i\ell}_{r}$ is assumed at time instant $se^{i\ell}_r$, for $r=1, \ldots, re^{i\ell}$.} The external arrivals represent block of arrivals of passengers for instance, due to the end of a football match in a place close to one of our stations. For the sake of readability and without loss of generality,  we will assume that $se^{i\ell}_0=0,\, se_{re^{i\ell}+1}^{i\ell}=\widehat{T}_\ell$, and $\Psi^{i\ell}_{0}=\Psi^{i\ell}_{re^{i\ell}+1}=0$, i.e., the first and last time instants of external arrivals coincide with the beginning and the end of the {extended} planning horizon, and the discontinuity jumps at those instants are null. Given a time instant $t \in [0,\widehat{T}_\ell]$, we use the following set of binary variables to determine {whether $t$ belongs to the interval} $[se^{i\ell}_r, se^{i\ell}_{r+1})$ for $ {r=0, \ldots, re^{i\ell}}$:
$$
\delta^E_{ri\ell}(t) = \left\{\begin{array}{cl} 1 & \mbox{if $t \in [se^{i\ell}_r, se^{i\ell}_{r+1})$,}\\
0 & \mbox{otherwise,}\end{array}\right. \quad i\in N_\ell,\, \ell\in L.
$$
Note that with these settings, the accumulated discontinuity flow jumps of external arrivals for a station $i\in N_\ell$ of line $\ell\in L$ can be modeled using the following constraints

\begin{equation}\label{DE}\tag{${\rm D}_{\rm E}$}
\begin{split}
 J^E_{i\ell} (t) = \dsum_{r=0}^{re^{i\ell}} \left(\dsum_{r^\prime\leq r} \Psi_{ir^\prime}^\ell\right) \delta^E_{r i\ell}(t),\quad  i\in N_\ell,\, \ell\in L,\\
 se^{i\ell}_r \delta^E_{ri\ell}(t) \leq t <  se^{i\ell}_{r+1} \delta^E_{ri\ell}(t)+\widehat{T}_\ell (1-\delta^E_{ri\ell}(t)),\quad r=0, \ldots, re^{i\ell},\, i\in N_\ell,\, \ell\in L,\\
\dsum_{r=0}^{re^{i\ell}} \delta^E_{ri\ell}(t)  =1,\quad  i\in N_\ell,\, \ell\in L,\\
\delta^E_{ri\ell}(t) \in \{0,1\}\quad r=0, \ldots, re^{i\ell},\, i\in N_\ell,\, \ell\in L.
\end{split}
\end{equation}

The reader may observe that the first constraint allows us to appropriately define $J^E$ by accumulating the flow jumps previous to a given instant $t$. The second set of constraints permits to determine the intervals in which $t$ lies and the third constraint ensures that only one of these intervals is identified for each $t$.

\subsection{Internal Arrivals\label{ss:IA}}

The main difference between internal and external block of arrivals is that in the latter, the breakpoints and the discontinuity flow jumps are known, while for internal arrivals, this information depends on the decision variables of the problem. {Recall that internal block of arrivals occur when passengers get off a train in an interchange station to transfer to another line. Thus, the time instants of those arrivals are part of the decision problem.} Therefore, we need to state a set of equations showing the existing {relationships} between the decision variables and the breakpoints and jumps in the flow of passengers they provoke.

In what follows we describe the modeling issues concerning the breakpoint  {times} and flow jumps of internal arrivals.

\begin{itemize}
\item {\bf Breakpoint {times} at line $\ell\in L$.}

Note that the instants{,  $si^{i\ell\ell^\prime}_r$,} in which an internal flow jump occurs by {the} arrivals of { trains coming from } line $\ell^\prime\in L$ to  {an} interchange station $i\in N_\ell \cap N_{\ell^\prime}$ {can be computed in terms of the time $t_i^{r\ell^\prime}$ as:}
$$
 si^{i\ell\ell^\prime}_r = t_{i}^{r\ell^\prime} - e_{i}^{r\ell^\prime}, \quad r=1, \ldots, \bar{k}_{\ell^\prime}
$$
{that is,} the {time} {instant} in which the $r$-th trip of line $\ell^\prime\in L$ departs to the next station ($ t_{i}^{r\ell^\prime}$) minus the waiting time at the station $i\in N_{\ell^\prime}$ ($e_{i}^{r\ell^\prime}$), that is, the instants in which the trains arrive to the interchange station of line $\ell'\in L$.

\item {\bf Discontinuity flow jumps at line $\ell\in L$ coming from line $\ell^\prime\in L$.}

The volume of the internal block of arrivals to an interchange station can  be derived from the passenger flows controlled by our decision variables. To compute it we will also need a set of values quantifying the proportion of passengers which want to change from one metro line to another one in all the interchange stations. Let $\tau_i^{\ell\ell^\prime}$ be the proportion of the passengers that get off a train in an interchange station $i\in N_\ell \cap N_{\ell^\prime}$  of the line  $\ell^\prime\in L$ to change to line $\ell \in L$. Thus, the flow jump at an instant $si^{i\ell\ell^\prime}_r$ for $r\in K_{\ell^\prime},\, \ell\in L, \, i \in N_\ell$ is given by
$$
\Phi_{ir}^{\ell\ell^\prime} =\left\{\begin{array}{ll}
 \tau_{i}^{\ell\ell^\prime} \dsum_{j<i} p_{ji}^{\ell^\prime}f_j^{r\ell^\prime}&\mbox{ if $(i,\ell^\prime)\in \overline{\mathcal{S}}$}\\
 \tau_{i}^{\ell\ell^\prime} \left(\dsum_{j<i} p_{ji}^{\ell^\prime}f_j^{r\ell^\prime}+\dsum_{j<i; j\in S_{\ell^\prime}} p_{ji}^{\ell^\prime}g_j^{r\ell^\prime}\right)&\mbox{ if $i\in S_{\ell^\prime},\, \ell^\prime\in LS$}\end{array}\right.
$$
\end{itemize}

Let $S^I_{i\ell\ell^\prime}=\left\{si^{i\ell\ell^\prime}_0, \cdots, si^{i\ell\ell^\prime}_{\bar{k}_{\ell^\prime}+1}\right\}$ be the set of sorted {breakpoint} {times} at {a given interchange} station $i\in N_\ell \cap N_{\ell^\prime}$ of line $\ell\in L$ caused by line $\ell^\prime\in L$. In order to make clearer the formulation, we will assume, w.l.o.g., that $si^{i\ell\ell^\prime}_0=0,\, si^{i\ell\ell^\prime}_{\bar{k}_{\ell^\prime}+1}=\widehat{T}_\ell$, and $\Phi_{i0}^{\ell\ell^\prime}=\Phi_{i\bar{k}_{\ell^\prime}+1}^{\ell\ell^\prime}=0$, i.e., the first and last time instants of internal arrivals coincide with the beginning and the end of the {extended} planning horizon, and the flow jumps at those times are null.

Now, we {are ready to model}  the internal flow jumps induced by line $\ell^\prime\in L$ onto the station $i\in N_\ell \cap N_{\ell^\prime}$ of the line $\ell$. {We proceed analogously} {as in} {the external arrivals} but incorporating a set of binary variables ($\delta^I$) to identify in which time interval between two consecutive breakpoints a given instant $t$ {belongs to}:
\begin{equation}\label{DI}\tag{${\rm D}_{\rm I}$}
\begin{split}
J^I_{i\ell\ell'}(t)  = \dsum_{r=0}^{\bar{k}_{\ell^\prime}} \left(\dsum_{r^\prime\leq r} \Phi_{ir^\prime}^{\ell\ell^\prime}\right) \delta^I_{r i\ell\ell^\prime}(t), \quad  i\in N_\ell \cap N_{\ell^\prime},\, \ell\in L, \ell^\prime\in L,\\
si^{i\ell\ell^\prime}_r \delta^I_{ri\ell\ell^\prime}(t) \leq t < si^{i\ell\ell^\prime}_{r+1} \delta^I_{ri\ell\ell^\prime}(t)+\widehat{T}_\ell(1-\delta^I_{ri\ell\ell^\prime}(t)),\; r=0, \ldots, \bar{k}_{\ell^\prime},\, i\in N_\ell \cap N_{\ell^\prime},\, \ell\in L, \ell^\prime\in L,\\
\dsum_{r=0}^{\bar{k}_{\ell^\prime}} \delta^I_{ri\ell\ell^\prime}(t)  =1, \quad  i\in N_\ell \cap N_{\ell^\prime,}\, \ell\in L, \ell^\prime\in L,\\
\delta^I_{ri\ell\ell^\prime}(t) \in \{0,1\}\; r=0, \ldots, \bar{k}_{\ell^\prime},\, i\in N_\ell \cap N_{\ell^\prime},\, \ell\in L, \ell^\prime\in L.
\end{split}
\end{equation}
{Note that}  the {two first sets of the above inequalities} are nonlinear since  $\mathbf{\Phi}${, $\mathbf{S^I}$}  and $\mathbf{\delta^I}$ are decision variables of our problem.

{As mentioned above}, the demand function $D_i^\ell(t)$ is only applied to a certain (finite) set of time {instants}. Then, {it} can be incorporated into the model by {adding} the variables $D_{i}^{k\ell} \equiv D_i^\ell (t_{i}^{k\ell})$,  $\delta^E$ and $\delta^I$, as well as  their corresponding linear and nonlinear constraints in \eqref{DE} and \eqref{DI}. We also include in the model \eqref{MINLP} { the following} two sets of valid inequalities:
 \begin{align}
    & \delta_{r i \ell}^{E}(t_i^{k+1\ell})\leq \sum_{r'=0}^{r}  \delta_{r^\prime i \ell}^{E} (t_i^{k\ell}),&  k \in K_\ell,\, ,r=0, \ldots, re^{i\ell},\,  i\in N_\ell,\, \ell\in L, \label{in:13}\\
     & \delta_{r i \ell\ell^\prime}^{I}(t_i^{k+1\ell})\leq \sum_{r'=0}^{r}  \delta_{r^\prime i \ell\ell^\prime}^{I} (t_i^{k\ell}),& k \in K_\ell, \; r=0, \ldots, \bar{k}_{\ell^\prime},\, i\in N_\ell \cap N_{\ell^\prime},\, \ell\in L, \ell^\prime\in L. \label{in:14}
    \end{align}
Inequalities \eqref{in:13} and \eqref{in:14}, ensure that trip $k+1$ {goes} through station $i$ after trip $k$. Although these inequalities are redundant with the rest of the constraints in the model, they considerably improve the performance of the solver over the branch-and-bound tree induced by the mathematical programming model.

\section{Case study: The Metrolab{\textregistered}  pilot experience\label{ex0}}
We illustrate the MINLP model introduced in Section \ref{sec:2} together with  the demand function described in Section \ref{sec:3} with a simplified version of the network of our pilot experience for Metrolab{\textregistered}.

Consider the simple network example depicted in Figure \ref{ex0:fig}. This network has a similar topology to the one provided by Metrolab{\textregistered} to calibrate our model in the pilot experience. It represents a small section of the Paris subway network.

\begin{figure}[h]
\begin{center}
\scalebox{0.8}{\begin{tikzpicture}[scale=0.8]
\coordinate(X1) at (-4,0);
\coordinate(X2) at (-1,0);
\coordinate(X3) at (0,0);
\coordinate(X4) at (2,0);
\coordinate(X5) at (5,0);

\coordinate(Y1) at (0,5);
\coordinate(Y2) at (0,3);
\coordinate(Y4) at (0,-2);
\coordinate(Y5) at (0,-4);

\node[circle,draw, inner sep=1pt](X-1) at (X1) {1};
\node[circle,draw, gray, inner sep=1pt](X-2) at (X2) {2};
\node[circle,draw, gray, inner sep=1pt](X-3) at (X3) {3};
\node[circle,draw, gray, inner sep=1pt](X-4) at (X4) {4};
\node[circle,draw, inner sep=1pt](X-5) at (X5) {5};

\node[circle,draw, inner sep=1pt](Y-1) at (Y1) {1'};
\node[circle,draw, gray, inner sep=1pt](Y-2) at (Y2) {2'};
\node[circle,draw, gray, inner sep=1pt](Y-4) at (Y4) {4'};
\node[circle,draw, inner sep=1pt](Y-5) at (Y5) {5'};

\draw (X-1)--node[above] {$3$} (X-2) ;
\draw[very thick, dashed] (X-2)--node[above] {$1$}(X-3);
\draw[very thick, dashed] (X-3)--node[above] {$2$}(X-4);
\draw (X-4)--node[above] {$3$}(X-5);

\draw (Y-1)--node[right] {$2$}(Y-2);
\draw[very thick, dashed] (Y-2)--node[right] {$3$}(X-3);
\draw[very thick, dashed] (X-3)--node[right] {$2$}(Y-4);
\draw (Y-4)--node[right] {$2$}(Y-5);
\end{tikzpicture}}
\end{center}
\caption{Network of Case study \ref{ex0}.\label{ex0:fig}}
\end{figure}
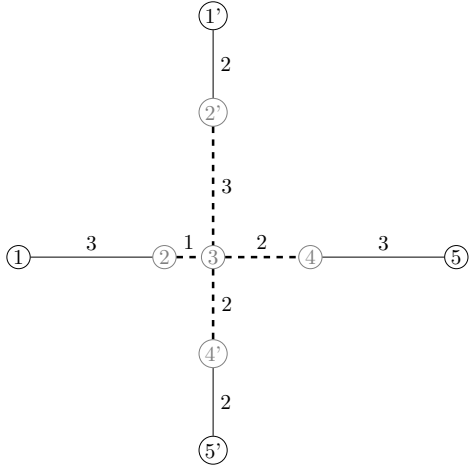

The network consists of {two \textit{bidirectional lines}}{(1-2-3-4-5 and 1'-2'-3-4'-5'}), each of them with five stations and sharing one of them (station $3$) which acts as an interchange station. Both lines allow short-turns on the planning. The edges of the lines which take part of the short-turn of each line are marked with dashed lines ($2-3-4$ for {one of the lines} and $2'-3-4'$ for {the other}). The travel time between consecutive stations is written next to each edge of the network. The stopping times are fixed to 30 seconds for stations different from the head and the final ones ($e_i^\ell$-parameters). Also the minimum safety time between consecutive trips is 2 minutes ($IS^\ell$-parameter). Hence, using our notation {$|L| = 4$} (two lines, but each of them in two directions), {$|LS|=4$ and $|LNS|=0$}. The lines will be numbered from 1 to 4, where $1$ and $2$ correspond with the horizontal line, left--right and right--left, respectively, and lines $3$ and $4$ are identified with the vertical line, up--down  and down--up, respectively. {The set of available capacities for the trains are $800$ and $1600$ passengers.}

For illustration purposes, we consider a planning horizon of $T=20$ minutes (from 7:30 to 7:50), and a maximum number of trips $K_\ell=7$ for the four lines.  We  assume that all the passengers that cannot get on a train because it is full await for the next train ($\alpha =1$) and that  the unitary penalty for persisting passengers is fixed to $\mu_1=0.1875$ (except for the last trip in which we fix it to $10 \times \mu_1=1.875$ to avoid a high excess of passengers at the end of the planning horizon). We also assume that $40\%$ of the passengers that get-off at station $3$ interchange to the other lines ($\tau$-parameter).

The {O-D matrix} and rewards per passenger are detailed in Tables \ref{table:exflows} and \ref{table:rewards}, respectively.

\renewcommand{\tabcolsep}{0.12cm}

\begin{table}[h]
\begin{center}
{\begin{tabular}{r|ccccc|r|ccccc|}
\multicolumn{1}{r}{$p_{ij}^\ell$} & {\tiny 1}     & {\tiny 2}     & {\tiny 3}     & {\tiny 4}     & \multicolumn{1}{c}{\tiny 5} & \multicolumn{1}{r}{$p_{ij}^\ell$} & {\tiny 1'}    & {\tiny 2'}    & {\tiny 3}     & {\tiny 4'}    & \multicolumn{1}{c}{\tiny 5'}\\
\cline{2-6}\cline{8-12}{\tiny 1}     & 0     & 0.40   & 0.35  & 0.20   & 0     & {\tiny 1'}    & 0     & 0.40   & 0.35  & 0.20   & 0\\
{\tiny 2}     & 0.40   & 0     & 0.60   & 0.35  & 0     & \hspace*{1cm} {\tiny 2'}    & 0.40   & 0     & 0.60   & 0.35  & 0 \\
{\tiny 3}     & 0.35  & 0.6   & 0     & 0.95  & 0     & {\tiny 3 }    & 0.35  & 0.60   & 0     & 0.95  & 0 \\
{\tiny 4}     & 0.20   & 0.35  & 0.95  & 0     & 1     & {\tiny 4'}    & 0.20   & 0.35  & 0.95  & 0     & 1 \\
{\tiny 5}     & 0.05  & 0.05  & 0.05  & 1     & 0     & {\tiny 5'}    & 0.05  & 0.05  & 0.05  & 1     & 0\\
\cline{2-6}\cline{8-12}\end{tabular}}
\end{center}
\caption{{O-D matrix} of Case study \ref{ex0}: Lines 1-2 (left) and 3-4 (right). \label{table:exflows}}
\end{table}

\renewcommand{\tabcolsep}{0.15cm}
\begin{table}[h]
\begin{center}
{\begin{tabular}{r|rrrrr|r|rrrrr|}
\multicolumn{1}{r}{$\gamma_{ij}^\ell$} & \multicolumn{1}{l}{{\tiny 1}} & \multicolumn{1}{l}{{\tiny 2}} & \multicolumn{1}{l}{{\tiny 3}} & \multicolumn{1}{l}{{\tiny 4}} & \multicolumn{1}{l}{{\tiny 5}} & \multicolumn{1}{r}{$\gamma_{ij}^\ell$} & \multicolumn{1}{l}{{\tiny 1'}} & \multicolumn{1}{l}{{\tiny 2'}} & \multicolumn{1}{l}{{\tiny 3}} & \multicolumn{1}{l}{{\tiny 4'}} & \multicolumn{1}{l}{{\tiny 5'}} \\
\cline{2-6}\cline{8-12}{\tiny 1} & 0     & 0.3   & 0.4   & 0.6   & 1     & {\tiny 1'} & 0     & 0.2   & 0.5   & 0.7   & 1\\
{\tiny 2} & 0.3   & 0     & 0.1   & 0.3   & 1     &\hspace*{1cm} {\tiny 2'} & 0.2   & 0     & 0.3   & 0.5   & 1 \\
{\tiny 3} & 0.5   & 0.2   & 0     & 0.2   & 1     & {\tiny 3}\, & 0.4   & 0.2   & 0     & 0.2   & 0 \\
{\tiny 4} & 0.6   & 0.3   & 0.1   & 0     & 0     & {\tiny 4'} & 0.7   & 0.5   & 0.3   & 0     & 0 \\
{\tiny 5} & 0.9   & 0.6   & 0.4   & 0.3   & 0     & {\tiny 5'} & 0.9   & 0.7   & 0.5   & 0.2   & 0 \\
\cline{2-6}\cline{8-12}\end{tabular}}
\end{center}
\caption{Rewards of Case study \ref{ex0}: Lines 1-2 (left) and 3-4 (right). \label{table:rewards}}
\end{table}

For the demand function, we consider the arrival rates of passengers per minute and the initial number of passengers at each station of each line detailed in Table \ref{table:demand}.

\begin{table}[h]
\begin{center}
{\begin{tabular}{|c|ccccc|ccccc|}
\cline{2-11}\multicolumn{1}{r|}{} & \multicolumn{10}{c|}{Lines}\\
\cline{2-11}\multicolumn{1}{r|}{} & \multicolumn{5}{c|}{$\ell=1$\phantom{'} }           & \multicolumn{5}{c|}{$\ell=2$}             \\
\hline
Stations {($i$)}      & 1\phantom{'}     & 2\phantom{'}      & 3     & 4\phantom{'}      & 5\phantom{'}      & 5\phantom{'}      & 4\phantom{'}      & 3     & 2\phantom{'}              & 1\phantom{'}       \\\hline
{$\beta_{0i}^\ell$} & 50    & 50    & 50    & 50    & 0\phantom{0}     & 50    & 50    & 50  & 50       & 0\phantom{0}     \\
{$\beta_i^\ell$}     & 10    & 100   & 120   & 90\phantom{0}    & 0     & 10    & 160   & 180 & 150   & 0     \\
\hline\multicolumn{1}{r|}{}     & \multicolumn{5}{c|}{$\ell=3$}           & \multicolumn{5}{c|}{$\ell=4$}   \\
\hline
Stations {($i$)}          & 1'  & 2'  & 3      & 4'      & 5'     & 5'      & 4'      & 3     & 2'   & 1'\\\hline
{$\beta_{0i}^\ell$}     & 50    & 50    & 50    & 50    & 50    & 50    & 50    & 50    & 50    & 50 \\
{$\beta_i^\ell$}          & 10    & 150   & 170   & 160   & 0     & 10    & 100   & 180   & 150   & 0 \\
\hline
\end{tabular}}
\end{center}
\caption{Coefficients of the Demand functions of Case study \ref{ex0}.\label{table:demand}}
\end{table}

The timetabling for the first line obtained after running our model in Gurobi 8.0 \cite{gurobi} is provided in Table \ref{exresults1}. The reported solution was obtained after 12 hours running time with a MIP GAP of $1.51\%$. There, we detail  for each trip ($k$, indicating with 'S' those trips which are short-turns ), its optimal capacity, the departing times of each train of each line (DepTime),  and also the flow estimations at each of the stages: the number of passengers that get off the train (Get-Off), the number of passengers that get-on the train ($f_i^{k\ell}$ for {whole trips} or $g_{i}^{k\ell}$ for short-turns), the excess of passengers ($h_{i}^{k\ell}$), the  passengers surplus of true trips ($x_i^{k\ell}$) {and the actual load of the train}.

\begin{table}[h]
{
\begin{center}
{\small\begin{tabular}{|c|c|c|c|c|c|c|c|}
\hline
 \textbf{$k$: Capacity} &$i$ & \textbf{DepTime} & \textbf{Get-Off} & $f_{i}^{k\ell}$ ($g_{i}^{k\ell}$) & $h_{i}^{k\ell}$ & $x_{i}^{k\ell}$& \textbf{Load} \\\hline
                       \multirow{5}{*}{1: 800} & 1     &  07:30:00 & 0.00  & 50.00 & 0.00  & 0.00  & 50.00 \\
\cline{2-8}             & 2     &  07:33:30  & 20.00 & 400.00 & 0.00  & 0.00  & 430.00 \\
\cline{2-8}             & 3     &  07:35:00  & 257.50 & 627.50 & 101.50 & 101.50 & 800.00 \\
\cline{2-8}             & 4     &  07:37:30  & 746.13 & 725.00 & 0.00  & 0.00  & 778.88 \\
\cline{2-8}             & 5     &  07:40:30  & 778.88 & 0.00  & 0.00  & 0.00  & 0.00 \\
\cline{1-8}       \multirow{3}{*}{2S: 1600} & 2     &  07:39:34 & 0.00  & 606.94 & 0.00  & 0.00  & 606.94 \\
\cline{2-8}             & 3     &  07:41:04 & 364.17 & 1231.59 & 0.00  & 0.00  & 1474.36 \\
\cline{2-8}             & 4     &  07:43:34 & 1474.36 & 0.00  & 638.18 & 0.00  & 91.93 \\
\cline{1-8}       \multirow{5}{*}{3: {\bf 0}} & 1     &  07:30:00 & 0.00  & 0.00  & 0.00  & 0.00  & 0.00 \\
\cline{2-8}             & 2     &  07:39:34 & 0.00  & 0.00  & 0.00  & 0.00  & 0.00 \\
\cline{2-8}             & 3     &  07:41:04 & 0.00  & 0.00  & 0.00  & 0.00  & 0.00 \\
\cline{2-8}             & 4     &  07:43:34 & 0.00  & 0.00  & 638.18 & 0.00  & 0.00 \\
\cline{2-8}             & 5     &  07:40:30 & 0.00  & 0.00  & 0.00  & 0.00  & 0.00 \\
\cline{1-8}       \multirow{3}{*}{4S: 800} & 2     &  07:43:12 & 0.00  & 364.02 & 0.00  & 0.00  & 364.02 \\
\cline{2-8}             & 3     &  07:44:42 & 218.41 & 603.05 & 0.00  & 0.00  & 748.66 \\
\cline{2-8}             & 4     &  07:47:12 & 748.66 & 0.00  & 1014.15 & 0.00 & 48.35 \\
\cline{1-8}       \multirow{5}{*}{5: {\bf 0}} & 1     &  07:30:00 & 0.00  & 0.00  & 0.00  & 0.00  & 0.00 \\
\cline{2-8}             & 2     &  07:43:12 & 0.00  & 0.00  & 0.00  & 0.00  & 0.00 \\
\cline{2-8}             & 3     &  07:44:42 & 0.00  & 0.00  & 0.00  & 0.00  & 0.00 \\
\cline{2-8}             & 4     &  07:47:12 & 0.00  & 0.00  & 1014.15 & 0.00 & 0.00 \\
\cline{2-8}             & 5     &  07:40:30 & 0.00  & 0.00  & 0.00  & 0.00  & 0.00 \\
\cline{1-8}       \multirow{5}{*}{6: 1600} & 1     &  07:46:15 & 0.00  & 162.60 & 0.00  & 0.00  & 162.60 \\
\cline{2-8}             & 2     &  07:49:45  & 65.04 & 655.05 & 0.00  & 0.00  & 752.61 \\
\cline{2-8}             & 3     &  07:51:15  & 449.94 & 1297.33 & 0.00  & 0.00  & 1600.00 \\
\cline{2-8}             & 4     &  07:53:45  & 1494.25 & 1494.25 & 109.44 & 109.44 & 1600.00 \\
\cline{2-8}             & 5     &  07:56:45  & 1600.00 & 0.00  & 0.00  & 0.00  & 0.00 \\
\cline{1-8}       \multirow{5}{*}{7: 800} & 1     &  07:50:00 & 0.00  & 37.40 & 0.00  & 0.00  & 37.40 \\
\cline{2-8}             & 2     &  07:53:30  & 14.96 & 373.99 & 0.00  & 0.00  & 396.43 \\
\cline{2-8}             & 3     &  07:55:00  & 237.48 & 641.05 & 0.00  & 0.00  & 800.00 \\
\cline{2-8}             & 4     &  07:57:30  & 747.38 & 446.03 & 0.00  & 0.00  & 498.65 \\
\cline{2-8}             & 5     &  08:00:30  & 498.65 & 0.00  & 0.00  & 0.00  & 0.00 \\
\hline
\end{tabular}}
\end{center}}
\caption{Results for {the first line}  of the Case study  \ref{ex0}. \label{exresults1}}
\end{table}
Note that trips $3$ and $5$ are fake trips (they have zero capacity), so the optimal number of trips for the planning is $5$. Two of them, the second and third true trips (trips $2$S and $4$S in Table \ref{exresults1}) are short-turns, while the remainder are {whole} trips with different capacities. We can  observe that fake trips  operate  on the line at the same time as the previous true trip. For instance, trip $3$ departs from station $1$, which doesn't belong to the short-turn, at the same time than the previous true whole trip (trip $1$) and from station $2$, which belongs to the short-turn, at the same time than the previous true trip for this station (short-turn trip $2$). Note also that the $\mu_1$-parameter considered for the last trip of the lines causes, in this case, that the excess of passengers is equal to $0$ at the end of the planning horizon.

The departing times of each true trip of all the lines from each of the stations are detailed in Figure \ref{fig:ex}. The horizontal axis is the time horizon, while the vertical axis represent spatial lengths (stations). We also report over each of the itineraries, their optimal capacities. Shortest itineraries running only over a subset of stations are short-turns. One can observe that the time difference between consecutive trips is not constant, as expected from the asymmetry of the lines and the transfer of passengers between them at different time instants. Then, as mentioned before, we deal with an aperiodic timetabling.

\begin{figure}[h]
\begin{center}
\includegraphics[scale=0.21]{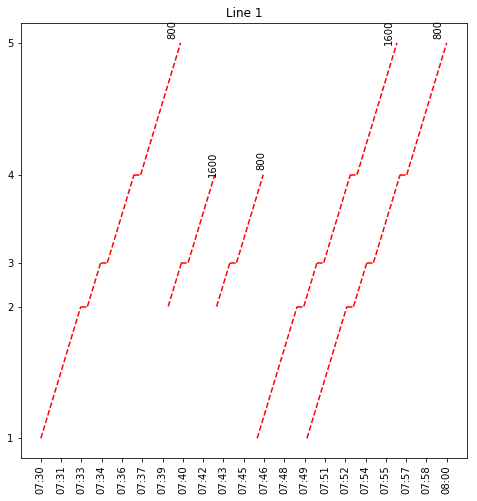}~\includegraphics[scale=0.21]{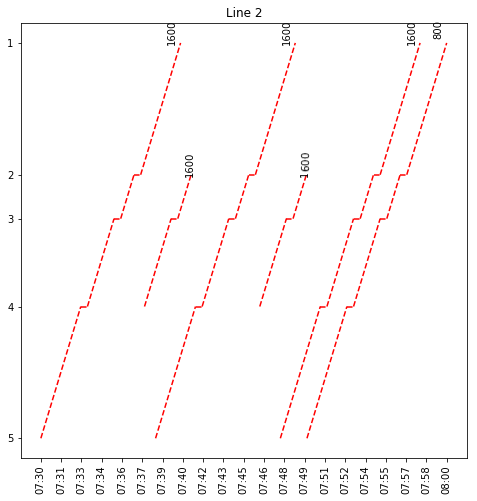}~\includegraphics[scale=0.21]{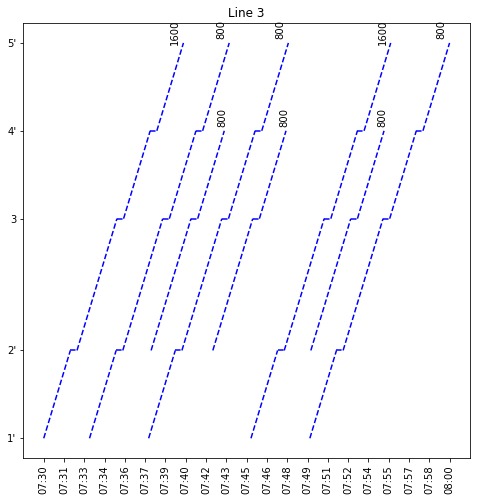}~\includegraphics[scale=0.21]{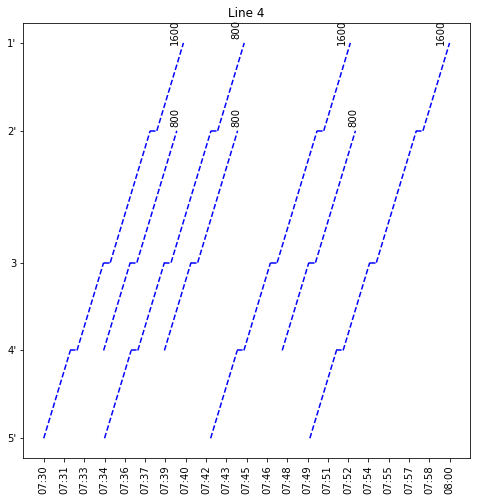}
\caption{{Timetables of the different trips for each line of the Case study of Section \ref{ex0}}.\label{fig:ex}}
\end{center}
\end{figure}

Furthermore, with the solution obtained, the estimated demands  at the departing times of each one of the trips at the interchange station are drawn in Figure \ref{fig:demands}. Observe that in that picture we draw in the horizontal axis, the departing times for such an interchange station at the different trips. In the vertical axis we depict the accumulated demand at those time instants. Note that the jumps in the demand induced by passengers that interchange to that line occur during the whole time interval between departing times of the interchange station, but we only account for  the accumulated amount when the new train {departs}.

\begin{figure}[H]
\begin{center}
\includegraphics[scale=0.2]{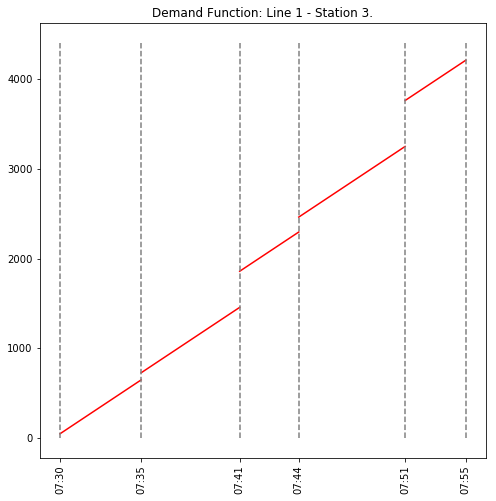}~\includegraphics[scale=0.2]{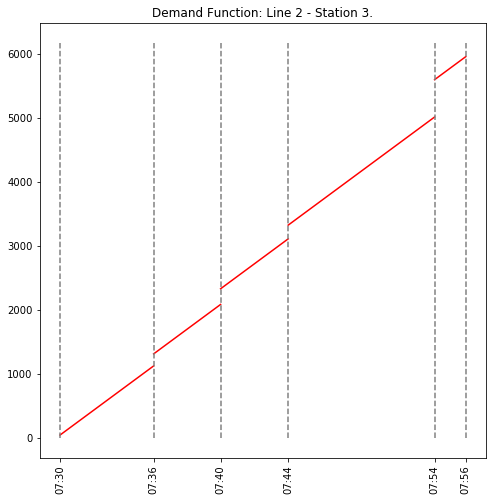}~
\includegraphics[scale=0.2]{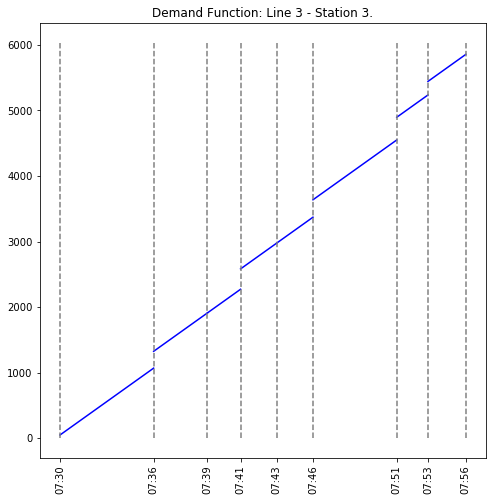}~\includegraphics[scale=0.2]{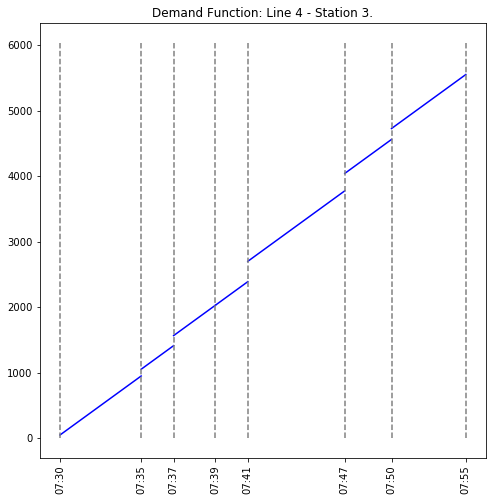}
\caption{Accumulated Demands obtained for the Case study of Section \ref{ex0}.\label{fig:demands}}
\end{center}
\end{figure}

A summary of the costs and rewards obtained for the reported solution is given in Table \ref{ex:summary}. Since an aggregated function of the costs is minimized in our model, the negative overall cost obtained ($-11404.88$) can be seen as a positive global reward for the optimal planning.

\begin{table}[!]
\begin{center}
\begin{tabular}{|c||cccc|c|}\hline
Costs/Rewards  & Line 1 & Line 2 & Line 3 & Line 4 & All lines \\\hline\hline
\eqref{obj:2} Complete Lines & 409.08 & 294.00 & 294.00& 294.00 & 1291.08\\
\eqref{obj:2} Short-Turns & 102.27 & 55.99 & 66.00 & 66.00 & 290.26\\
\eqref{obj:4a} {Complete Lines }&  1736.16 & 2514.53&3365.67&2814.05 & 100430.41\\
\eqref{obj:4a} { Short-Turns }& 541.70 &534.52 & 743.69& 780.74 & 2448.67\\
\eqref{obj:5} &  44.85 & 0 & 0& 0 & 44.85\\\hline
OVERALL COST & -1721.65& -2699.06 & -3749.37 &  -3234.79 &  {\bf -11404.88}\\\hline
\end{tabular}
\end{center}
\caption{Summary of solution obtained for the Case study of Section \ref{ex0}.\label{ex:summary}}
\end{table}

\bigskip

Problem \eqref{MINLP} depends on the demand that occur in each line and {on the correlation between lines induced by interchanging stations. For that reason, to have an exact model the demands by lines must be considered as variables and hence, the problem that considers all the lines at the same time becomes very hard by obvious reasons. Among them, the problem is a Mixed Integer Non Linear Programming Program whose continuous relaxation is neither convex nor concave because of the bilinear constraints in \eqref{DE} and \eqref{DI} (which are needed to be linearized, turning into weak continuous relaxations of the problem); and in addition, the number of variables increases considerably with each line that is jointly considered in the problem. Therefore, in order to provide an efficient procedure which is able to handle the model for realistic instance sizes in reasonable time (recall that the pilot Case study  \ref{ex0}} required 12 hours of CPU time to obtain a nearly optimal solution), we propose  a Math-Heuristic approach. This algorithm solves sequentially Problem \eqref{MINLP} for one line, say $\ell\in L$, at a time  (by fixing appropriately the demand parameters after each run) to approximate the actual solution of the global problem. We will denote such a problem as Problem \eqref{MINLP}$_\ell$.  This Math-Heuristic approach is developed in the following section.

\section{A Math-Heuristic Approach}
\label{sec:4}

\begin{figure}[h]
\begin{center}
\scalebox{0.5}{\tikzstyle{startstop} = [rectangle, rounded corners, minimum width=3cm, minimum height=1cm,text centered, draw=black]
\tikzstyle{elli} = [ellipse, minimum width=2.25cm, minimum height=1cm,text centered, draw=black]

\tikzstyle{init} = [rectangle, minimum width=3cm, minimum height=1cm,text centered, draw=black]

\tikzstyle{decision} = [diamond, node distance=3cm,
     minimum height=2em,text width=3.8cm, text badly centered,  draw=black, fill=gray!10]
\tikzstyle{arrow} = [thick,->,>=stealth]
\tikzstyle{arrow1} = [thick,<-,>=stealth]
\tikzstyle{line} = [draw, -latex']

\begin{tikzpicture}[node distance=2.25cm]
\node (start1) [init] {{\bf Initialization:} $\widehat{L}\leftarrow L$, $L_0\leftarrow\widehat{L}$};
\node (start1a) [startstop, below of=start1,minimum width=1.6cm] {$\ell \in L_{it}$};
\node (start2) [startstop, below of=start1a,fill=gray!10] {Update Demands};
\node (dec1) [decision,aspect=2, below of=start2, node distance=3cm] {Solve \eqref{MINLP}$_\ell$};
\node (dec1a)  [startstop, below of=dec1,minimum width=2cm, node distance=1.8cm]  {$L_{it} \leftarrow L_{it}\backslash\{\ell\}$};
\node (start3) [elli, below of=dec1a,node distance=2.25cm] {Stabilization{\large\bf ?}};
\node (start4) [startstop, below of=start3] {$\widehat{L}\leftarrow\widehat{L}\backslash\{\ell\}$};
\node (start4a) [elli, right of=start3, node distance =4cm] {$L_{it}\neq\emptyset${\large\bf ?}};
\node (start5) [elli, right of=start4a, node distance =5cm] {$it < \texttt{maxit}$ \& $\widehat{L}\neq\emptyset${\large\bf ?}};
\node (dec2) [decision,aspect=2, below of=start5, node distance=5cm] {{\bf Output:}\\
Solve \eqref{MINLP} for fixed binaries.};

\draw[arrow] (start1)-- node[right]{$it\leftarrow0$} (start1a);
\draw[arrow] (start1a)--  (start2);
\draw[arrow] (start2)-- (dec1);
\draw[arrow] (dec1a)-- (start3);
\draw[arrow] (start3)-- node[right] {yes} (start4);
\draw[arrow] (start3)-- node[above] {no} (start4a);
\draw[arrow] (start4)--(start4a);
\draw[arrow] (start5)-- node[right] {no} (dec2);
\draw[arrow] (start4a) -- node[right] {yes} ($(start2.east)+(2.55cm,0)$) -- node[below] {$\ell \in L_{it}$}   (start2);
\draw[arrow] (start4a)-- node[above] {no}  (start5);
\draw[arrow] (start5) -- node[right] {yes} ($(start1a.east)+(8.3cm,0)$) -- node[below] {$it\leftarrow it+1$}  node[above] {$L_{it}\leftarrow\widehat{L}$}  (start1a);
\end{tikzpicture}}
\end{center}
\caption{Flowchart of the Math-Heuristic approach.\label{fc}}.
\end{figure}
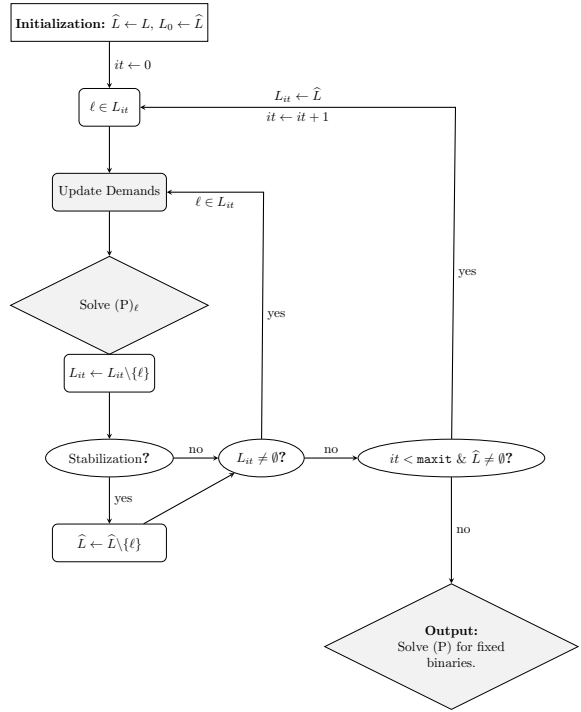

In this section we propose an optimization-based heuristic approach to solve \eqref{MINLP}. Note that the MINLP formulation stated for the problem involves many continuous and discrete variables, as well as a large number of constraints. In particular, the number of $\delta$-variables that allows us to model the internal jumps of the demand function, described in the previous section, {is} $O(K^2  N^{Ic} L^2)$, where $K$ is the maximum number of trips and $N^{Ic}$ is the number of interchange stations. That bound implies a high complexity to solve even small-size instances of the problem.

In our algorithm we apply a \textit{divide-and-conquer} strategy by splitting the models by lines, that is, we solve sequentially the model for a single line $\ell\in L$ at a time. The output of that solution together with the solutions of those previously solved are passed to the next line to be solved. Hence, the flow information provided by the lines in $L \backslash \{\ell\}$ (transfer passenger flows and train arrival times at the interchange stations) are considered as input for  solving the problem \eqref{MINLP}$_\ell$ for the fixed line $\ell \in L$. In this way, the demand function is updated and a new line is solved using the same approach. This procedure continues until a stopping rule evaluating the progress of the overall costs is satisfied, i.e. when the solution stabilizes or until a maximum number of iterations is reached.  This approach allows us to avoid considering all the $\delta$-variables together in a single model, which is one of the main drawbacks of the exact algorithm. A compact flow chart of the algorithm is drawn in Figure \ref{fc}.

In the initial iteration we consider the whole set of lines to solve, $\widehat{L}=L$, and the  time instants in which the demand is affected by passengers coming from an interchange station, $\mathbf{S^I}$ (it is assumed that, for the first line solved, there are no passengers {coming} from other lines). Then, for each line $\ell\in L$, we update {the demand function taking into account} the flows and breakdown times obtained when solving the previous lines and then, we solve the problem for the single line $\ell$ with such an input information. Observe that after the initial {line is solved} (in which the connecting flows are set to zero), the optimal flows give us information about the passengers that get on and {get off} from the trains of the line, what affect the rest of the lines.

Once all the lines have been solved, we repeat the procedure $\texttt{maxit}$ times, to avoid getting trapped in locally optimal solutions. Let $z^{\ell(it)}$ be the objective function value of problem \eqref{MINLP}$_\ell$  at iteration $it$.    If at an iteration $it <\texttt{maxit}$, this value stabilizes, that is, the relative deviation in the objective value with respect to the previous iteration of line $\ell\in L$, namely $\frac{|z^{\ell(it)}-z^{\ell(it-1)}|}{z^{\ell(it)}}$ is not significant enough (smaller than a tolerance $\varepsilon$), the corresponding solution for this line is kept fixed for the remainder iterations. Thus, the iterative procedure terminates either, after the iteration when the solution stabilizes with some degree of accuracy, $\varepsilon$,  (for all the lines) or when a maximum number of iterations is reached. Finally, to obtain a feasible solution (upper bound), we solve \eqref{MINLP} by fixing the binary variables obtained during the heuristic procedure. Such a problem becomes a continuous Linear Programming problem which is efficiently solved.

\section{Experiments \label{sec:5}}

\begin{figure}[h]
\input{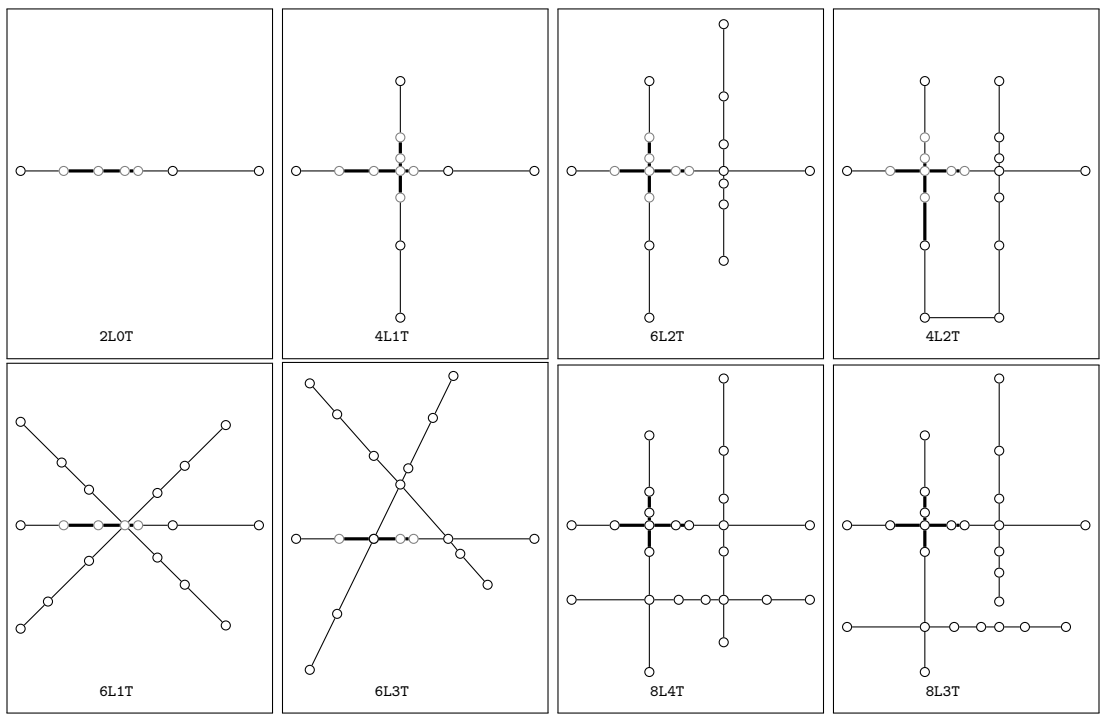}
\caption{Network Topologies for the Experiments.\label{fig:all}}
\end{figure}

We have run a series of computational experiments in order to test the model and the performance of our {two} approaches. We have considered the networks topologies drawn in Figure \ref{fig:all}. The names of the networks are in the form $X$\texttt{L}$Y$\texttt{T} where $X$ is the number of lines and $Y$ the number of interchange stations. We consider networks with $2, 4, 6$ and $8$ lines, each of them with $7$ stations except \texttt{4L2T} that contains 14 stations. For each of the $8$ networks, we tested our approach in two situations: no short-turns are allowed (that is, with our notation $|LS|=0$ and short-turns are allowed in some of the lines (thickest lines/gray nodes correspond with allowed short-turns for the lines) to evaluate the effect of short-turning.  The construction of the different networks is motivated by the representation of the most common situations in real-world subway networks {(see \cite{lap11})}. {The description of the networks is provided in Table \ref{table:lines}.

\begin{table}[h]
{\small\begin{center}
\begin{tabular}{|c|p{11.5cm}|}\hline
Network & Description\\\hline
\texttt{2L0T} & One bidirectional line with no interchange stations.\\\hline
\texttt{4L1T} &Two bidirectional lines with a single interchange station.\\\hline
\texttt{6L2T} & Three bidirectional lines. One of the bidirectional lines with two interchange stations, each of them interchanging for a different  bidirectional line (2-by-2).\\\hline
\texttt{4L2T} &{Two bidirectional lines.} One of the bidirectional lines with two interchange stations, each of them interchanging from the same second bidirectional line.\\\hline
\texttt{6L1T} &Three bidirectional lines with one interchange station for different lines (3-by-3).\\\hline
\texttt{6L3T} &Three bidirectional lines with three 2-by-2 interchange stations.\\\hline
\texttt{8L4T} &Four bidirectional lines with four 2-by-2 interchange stations.\\\hline
\texttt{8L3T} &Four bidirectional lines with three 2-by-2 interchange stations.\\\hline
\end{tabular}
\end{center}}
\caption{Description of the network topologies considered in our computational experiments.\label{table:lines}.}
\end{table}

Our set of networks include the one proposed by Metrolab{\textregistered} (\texttt{4L1T}) to analyze the viability of automatizing the management of the Paris subway network. Actually, the parameters of these networks were designed based on that case and are available at \url{http://bit.ly/InputData_SubwayPlanning}. All the models were coded in \texttt{Python 3.6}, and solved using \texttt{Gurobi 8.0}~\cite{gurobi} in a Mac OSX with an Intel Core i7 processor at 3.3 GHz and 16GB of RAM.

\subsection{Results}
\begin{table}
\begin{center}
{\small\begin{tabular}{|c|c|cc|cc|c|}\cline{3-6}
\multicolumn{ 2}{c}{} & \multicolumn{ 2}{|c}{{\bf Math-Heuristic}} & \multicolumn{ 2}{|c|}{{\bf MINLP} \eqref{MINLP}} & \multicolumn{1}{c}{} \\
\hline
{\bf Network} &   $|{\bf LS}|$      & {\bf BestObj} & {\bf CPU} & {\bf BestObj} & {\bf CPU} &     {\bf GAP} (\%)         \\
\hline\hline
    \multirow{2}{*}{\texttt{2L0T}} &          0 &     145702 &      $< 0.1$ &     145573 &          7 &       0.09 \\

           &          2 &     114145 &         11 &     112916 &      \texttt{TL}
           &       1.08 \\\hline
    \multirow{2}{*}{\texttt{4L1T}} &          0 &     206729 &        132 &     206242 &       \texttt{TL}
     &       0.24 \\
           &          4 &     152890 &        631 &     152016 &       \texttt{TL}
           &       0.57 \\\hline
    \multirow{2}{*}{\texttt{4L2T}} &          0 &     348267 &       3522 &     347224 &       \texttt{TL}
    &       0.30 \\
           &          4 &     333102 &       4892 &     332665 &       \texttt{TL}
           &       0.13 \\\hline
    \multirow{2}{*}{\texttt{6L1T}} &          0 &     276961 &       1130 &     276545 &       \texttt{TL}
    &       0.15 \\
           &          2 &     235488 &       1521 &     234589 &       \texttt{TL}
           &       0.38 \\ \hline
    \multirow{2}{*}{\texttt{6L2T}} &          0 &     249038 &        606 &     248854 &       \texttt{TL}
   &       0.07 \\
           &          4 &     203080 &       2882 &     203080  &\texttt{TL}
           &       0.00 \\\hline
    \multirow{2}{*}{\texttt{6L3T}} &          0 &     248988 &        520 &     248979 &       \texttt{TL}
    &     $< 0.01$ \\
           &          2 &     217004 &       2688 &     216909 &       \texttt{TL}
           &       0.04 \\\hline
    \multirow{2}{*}{\texttt{8L3T}} &          0 &     404627 &        432 &     404147 &       \texttt{TL}
   &       0.12 \\
           &          4 &     362916 &       1467 &     362908 &       \texttt{TL}
           &     $< 0.01$ \\\hline
    \multirow{2}{*}{\texttt{8L4T}} &          0 &     404469 &       1191 &     404211 &       \texttt{TL}
    &       0.06 \\
           &          4 &     374067 &       1144 &     374064 &       \texttt{TL}
           &     $< 0.01$ \\\hline
\end{tabular}}
\caption{Computational Results \label{Table:results}}
\end{center}
\end{table}

The results of our computational experience are shown  in Table \ref{Table:results}. In this table the results are organized in  {five} blocks of columns. The first block reports the network topology (using the name convention explained above).
The second block , ``$|{\bf LS}|$'', stands for the number of allowed short-turns on the considered network.
The {third} block  gathers the results obtained with the Math-Heuristic approach described in Section \ref{sec:4}, whereas the {fourth} one reports the corresponding results obtained with the exact MINLP Formulation  described in Section \ref{sec:2}. In both cases, `{`\textbf{BestObj}''} is the value of the objective function and  ``\textbf{CPU}'' is the time, in seconds, required to meet the stopping criterion or to reach the time limit (a maximum CPU time of 10 hours has been set  for solving the MINLP formulation), respectively. {Those instances for which the time limit was reached are indicated with \texttt{TL} in that column, in whose case, the best obtained solution is reported.} Finally, the block
``\textbf{Gap (\%)}'' is the percentage gap between  the objective function value obtained in the Math-Heuristic approach and the objective function value obtained when solving the MINLP Formulation.
It is worth noting that to improve the performance of the  MINLP solver, we have initialized it with an initial solution given by the one provided by the  Math-Heuristic approach. This implies that the solution of the exact approach is always as good as the Math-Heuristic and the \%Gap is a measure of the improvement given by the exact method over the proposed heuristic.

The results reported in Table \ref{Table:results} show the remarkable performance of our Math-Heuristic algorithm. In all cases, it achieves rather good solutions with gaps with respect to the exact method smaller than $1.08\%$ and in rather competitive computing times (one order of magnitude smaller than the exact MINLP solver). Note that the time limit of 10 hours was reached in all the instances except in the simplest one. The results are particularly exceptional in the most complex topologies using  eight lines where the improvement of the exact method is negligible: gaps are smaller than 0.12\% and computing times of the Math-Heuristic are, on average, around 2\% of those required by the MINLP. Based on our experiment we conclude that the Math-Heuristic algorithm is a good compromise to solve the line planning and timetabling problem considered in this paper. Apart from the relative deviations between the solutions obtained with the math-heuristic and the exact approach, Gurobi was not even able to solve the problems within the time limit in all the instances except the simplest one, while our math-heuristic is able to obtain good quality feasible solution of the problem in reasonable CPU times. Finally, it is worth mentioning that the consideration of short-turns in the network considerably increases the computational effort needed to solve the problem, but at the same time reduces the overall cost considered in our objective function (an average reduction of $15\%$ was observed under the considered parameters). For instance, the exact approach required 7 seconds to solve to optimality the instance \texttt{2L0T}, without short-turns. However, it was not able to certify optimality for the same instance with short-turns, within the time limit of 10 hours.  The comparison with respect to the Math-Heuristic is less dramatic and we report an average relative deviation of $57\%$ in CPU time  between the instances with and without short-turnings.

\section{Conclusions and Future work}
\label{sec:concl}
This paper considers a general  model for line planning and timetabling problems on subway networks that was originated by a collaboration with Metrolab{\textregistered}, a French R\&D company analyzing the viability of automating the metro of Paris. We propose a new Mathematical Programming- based decision making tool in this context. A number of different elements have been taken into account in the construction of both, objective function and constraint set of the proposed integrated line planning and timetabling model. The model aims to reflect as most real factors as possible of this complex transportation environment, but still being possible its computational suitability. Our outcome is a  Mixed Integer Non Linear Programming formulation  taking into account several quality and operation measures for feasible planing. We incorporate to the model an integrated cost- and passenger-oriented objective function, time-dependent demands and interchange stations together with several strategic decisions on the planning, as the starting times for each trip, the activation of short-turns, the determination of the number of trips in a journey, the selection of capacities for the trips, and the determination of the optimal number of trips. This representation of the problem may allow us to use the optimization model as a (\emph{what-if}) tool to assess potential technological or managerial innovations that may be introduced into the system. In order to assess such innovations, one can decouple the whole transportation system and consider just a (small) section of it where the main effects of such a new element produce a higher impact. Such a reduction process resembles the one followed in our paper to illustrate the case study. Apart from using a MINLP solver for solving small instances of the problem, we also develop a novel Math-Heuristic algorithm which allows us to solve realistic instances. Extensive computational experiments on different network topologies based on the one provided by Metrolab{\textregistered} on the metro of Paris, show that the Math-Heuristic algorithm performs remarkably well as compared to the exact resolution of the Mathematical Programming formulation, being an adequate tool to solve this type of problems.

Many future enhancements of the contributions made here could be cited now. Some of them could refer to the integration of new elements to the model such as variable speed for the trains or more general demand functions to model the number of passengers awaiting at the stations. Also, the stochastic nature of the parameters of the problem could be managed by incorporating the uncertainty on the passenger flows to the problem, and then using tools from Stochastic Programming to derive a model and different solution approaches for the problem. On the other hand, an interesting research line could be to obtain theoretical conditions under which one was able to ensure the convergence of our iterative numerical approach. It would be interesting to find conditions ensuring the convergence of the Math-Heuristic approach, even in simpler models, since they may give rise to new versions of the algorithm integrating different stopping rules or ways to decouple the system in parts.


\bibliographystyle{informs2014trsc}

\end{document}